\documentclass[12pt,twoside]{article}
\usepackage[latin1]{inputenc}
\usepackage[english]{babel}
\usepackage{amsmath,amsfonts,amssymb,amsthm}
\usepackage{graphicx}
\usepackage{epsfig}
\usepackage{caption}
\usepackage{float}
\usepackage{hyperref}
\usepackage{soul}
\usepackage{enumitem}
\allowdisplaybreaks

\newtheorem{theorem}{Theorem}[section]

\newtheorem{proposition}{Proposition}[section]
\newtheorem{corolario}{Corollary}[section]
\newtheorem{lema}{Lemma}[section]
\newtheorem{remark}{Remark}[section]

\newcommand{\N}{\mathbb{N}}
\newcommand{\V}{\text{Var}}

\newcommand{\tZ}{\tilde{Z}}
\newcommand{\tT}{\tilde{T}}
\newcommand{\tF}{\tilde{F}}
\newcommand{\tM}{\tilde{M}}
\newcommand{\ts}{\tilde{s}}
\newcommand{\tp}{\tilde{p}}
\newcommand{\tm}{\tilde{m}}
\newcommand{\trho}{\tilde{\rho}}
\newcommand{\talpha}{\tilde{\alpha}}
\newcommand{\tx}{\tilde{x}}

\newcommand{\tphi}{\tilde{\phi}}

\newcommand{\cuad}{\begin{flushright}\vspace{-2ex}$\Box$\vspace{-2ex}\end{flushright}}

\newenvironment{Prf}[1][\unskip]{%
\par
\noindent
{\textbf{Proof of #1}}\newline
\vspace{-2ex}\noindent{}\newline}\cuad

\usepackage[left=2cm,right=2cm,top=2cm,bottom=2cm]{geometry}
\author{Cristina Guti\'errez\footnote{Department of Mathematics, University of Extremadura, C\'aceres, Spain. E-mail address: \url{cgutierrez@unex.es}. ORCID: 0000-0003-1348-748X.} \and Carmen Minuesa\footnote{Department of Mathematics, Autonomous University of Madrid, Madrid, Spain. E-mail address: \url{carmen.minuesa@uam.es}. ORCID: 0000-0002-8858-3145.}}

\title{A predator-prey two-sex branching process}

\begin{document}

\maketitle

\begin{abstract}
In this paper, we introduce a two-sex controlled branching model to describe the interaction of predator and prey populations with sexual reproduction. This process is a two-type branching process, where the first type corresponds to the predator population and the second one to the prey population. While each population is described via a two-sex branching model, the interaction and survival of both groups is modelled through control functions depending on the current number of individuals of each type in the ecosystem. We provide necessary and sufficient conditions for the ultimate extinction of both species, the fixation of one of the species and the coexistence of both of them. Moreover, the description of the present predator-prey two-sex branching process on the fixation events can be performed in terms of the behaviour of a one-type two-sex branching process with a random control on the number of individuals, which is also introduced and analysed.
\end{abstract}



\section{Introduction}\label{sec:Introduction}

Predator-prey models have been widely studied in the literature since the introduction of the first  model in \cite{Lotka}. This process aimed at modelling the trophic interactions between both species and as later in \cite{Volterra}, non-linear differential equations were used to describe the interaction of a predator-prey dynamic system. Since then, many modifications and new models have been defined trying to adapt the peculiarities observed in the real world as faithfully as possible.  Some of the most recent papers on this topic are, for example, \cite{Lehtinen2019}, \cite{Lois2020}, \cite{Gao2020}, \cite{Yan2020},  and \cite{Sasmal2020}, which deal with theoretical questions, or  \cite{Fay2006} and  \cite{Grabowsky2019}, which focus on real world settings. 

All those papers use deterministic models based on ordinary differential equations (ODEs) to model the predator-prey dynamics. However, the interaction between both species can be seen as a stochastic system. In this sense, \cite{Hu2019} conducts an interesting study about the dynamics of deterministic and stochastic models for a predator-prey system, where the predator species suffers from a parasitic infection. The deterministic model is an ODE model while the stochastic model is derived by means of continuous time Markov chains. The theory of branching processes is used to estimate the probability of a disease outbreak and the probability of the prey species invasion.

In the context of branching processes several publications have tackled predator-prey system modelling. We stand out the pioneer work of  \cite{Hitchcock-1986} where two predator-prey models in continuous time are considered: a host-parasite system and a predator-prey process for which the predator birth rate is not directly associated with a prey death. A further study of the process in \cite{Hitchcock-1986} is performed in \cite{Ridler-Rowe-1988}. In \cite{Coffey-Buhler-1991}, a two-type branching process in discrete time is introduced to describe two populations, predators and preys, living in the same ecosystem. In the model, the evolution of the predator population is independent of the  population size of preys, and the number of preys at each generation is given by the number of prey offspring minus the number of preys that have been captured and killed by the predators. The analysis of this process continues in \cite{Alsmeyer-1993} and it focusses on the study of necessary and sufficient conditions for the fixation of both populations. Later, \cite{Coffey-1995} analyses a counterpart of the predator-prey model in \cite{Coffey-Buhler-1991} in the continuous-time setting by making use of linear birth-death processes. More recently, in \cite{Chapron2008} a branching model is used to analyse and compare the influence of the habitat lost, poaching or drop of preys in tiger populations. Moreover, it is important to mention the work of \cite{Durrett2010} where a branching random walk is used to describe predator-prey populations. 

However, none of the previous models takes into account the fact that many animal species have a sexual reproduction. This fact plays a key role in the evolution of the species because their development not only depends on the number of individuals of the species, but also on the total number of females and males and the type of mating of the populations. Many real cases have been reported on the interaction between preys and predators where both populations are formed by females and males which mate and procreate by means of sexual reproduction. For example, \cite{Morrison2017} studied populations of sea lions and penguins as predators and preys, respectively, and in \cite{Muller2019}, a social network of giraffe populations is studied bearing in mind the presence of lions which preferentially prey on giraffe calves.

Our aim in this paper is to introduce a predator-prey stochastic process by using the two-sex branching processes theory (see \cite{daleyA}) in order to fill this niche in the literature on predator-prey models. Therefore, the novelty of the present paper is to consider a two-type two-sex discrete-time branching process in order to model the predator-prey interaction of populations of females and males with sexual reproduction.  Moreover, we focus on situations where each female mates with only one male, whenever there are males in the population, but a same male could mate with more than one female. This type of mating is called promiscuous mating and we can observe it in many examples in the nature (see, for example, \cite{Norris1988}, \cite{Furtbauer2011}, \cite{Balme2012}, \cite{Lifjeld2019}, \cite{Wightman2019} or \cite{Lee2019}).

The two-sex branching process allows us to describe the generation-by-generation evolution of the number of individuals of a predator-prey system in certain environment. The definition of the model is based on the following assumptions: in each generation, females and males of each species mate and form couples by promiscuous mating, and each of those couples is assumed to give birth to some number of individuals. Nevertheless, the survival of all these individuals to form couples and reproduce is constrained due to the interaction between both species. Thus, some preys could be captured and killed by the predators to feed themselves and some predators could die due to lack of food supplies. As a result, the couples of each species of the following generation will be formed from the females and males that have survived. 

The second aim of this work is to study how the number of individuals of each species evolves over successive generations. We examine conditions for one of the species -predator or prey- to become extinct or to have a positive probability of survival. Moreover, we are interested in studying sufficient  conditions for the coexistence of both species assuming two different scenarios. First, we analyse the destiny of the species assuming that a predator population with limitless appetite is introduced and regarded as an invasive species in a geographically isolated area where its only food supply is the prey population. Second, we assume that predators with limited appetite could survive without any prey due to the presence of other food resources. This consideration on the number of preys consumed per predator will give the preys a chance to survive although they might coexist with a large predator population. 

This paper is divided into 7 sections. Apart from this introduction, a controlled two-sex branching process is introduced in Section \ref{sec:C2SBP}. This model will be useful to facilitate the understanding of the predator-prey two-sex branching model which is presented in Section \ref{sec:Definition}. The basic properties and main results about the conditional moments of the variables involved in the predator-prey two-sex branching model are studied in Section \ref{sec:BasicProperties}. In Sections \ref{sec:Isolated} and \ref{sec:Non-isolated}, we provide conditions for the ultimate extinction of the entire predator-prey system, for the coexistence of both species or for the fixation of one of them; in Section \ref{sec:Isolated} we focus on the case when both species live together in an isolated environment, whereas in Section \ref{sec:Non-isolated} we consider that they might coexist with other animal species in a non-isolated ecosystem. To conclude, a discussion about the main results obtained is reported in Section \ref{sec:Discussion}. In order to facilitate the reading, the proofs of the results are gathered in some appendices.

\section{A controlled two-sex branching process}\label{sec:C2SBP}

Before providing the definition of the predator-prey two-sex branching process, we begin by introducing a simpler branching process that will play an important role in the analysis developed in the following sections. The model is a generalization of the two-sex branching process introduced by \cite{daleyA} with the novelty that not all individuals generated by the couples of the previous generation participate in the mating process. Only some selected individuals -males and females- are chosen to formed couples. The selection of those individuals is done through random control functions. Thus, this model represents a combination of the two-sex branching process and the controlled branching process with random control functions on the total number of individuals. First, we provide the probabilistic definition of the model.

Let us consider two independent families of random variables (r.v.s) $\{\xi_{ni}:n\in\N_0,i\in\N\}$ and $\{\varphi_n(k):n,k\in\N_0\}$. The former is assumed to be a family of independent and identically distributed (i.i.d.) r.v.s. The latter is assumed to be a family of r.v.s such that $\{\varphi_n(k):k\in\N_0\}$, $n\in\N_0$, are independent stochastic processes with the same one-dimensional distribution. Let us also consider the sequences $\{Y_n\}_{n\in\N_0}$ and $\{X_{n}\}_{n\in\N}$ defined as:
\begin{equation}\label{def:C2SBP}
Y_0=y_0 \in\N,\quad X_{n+1}=\sum_{i=1}^{Y_n}\xi_{ni},\quad Y_{n+1}=L(\bar{F}_{n+1},\bar{M}_{n+1})\qquad n\in\N_0,
\end{equation}
where $\N_0=\N\cup\{0\}$, $L:\N_0^2\rightarrow \N_0$ is a  deterministic mating function and $(\bar{F}_{n+1},\bar{M}_{n+1})$ is a random vector that follows a multinomial distribution with parameters $y$ and $(\lambda,1-\lambda)$ conditionally on $\{\varphi_{n+1}(X_{n+1})=y\}$, with $0<\lambda<1$. The process $\{Y_n\}_{n\in\N_0}$ is known as \emph{two-sex branching process with random control on the total number of individuals (BBPCI)}. Moreover, the empty sum in \eqref{def:C2SBP} is assumed to be 0 and the mating function is assumed to be monotonic and non-decreasing in each argument and it satisfies the following conditions:
\begin{enumerate}[label=(A\arabic*),ref=(A\arabic*)]
\item $L(x,y)\leq x$, for each $x,y\in\N_0$.\label{cond:A2}
\item $L(0,y)=L(x,0)=0$, for each $x,y\in\N_0$.\label{cond:A3}
\end{enumerate}

\vspace{2ex}

Intuitively, the variable $Y_n$ represents the number of couples at generation $n$, while the variable $X_{n}$ denotes the total number of individuals -females and males- at generation $n$. The functions $\varphi_n(\cdot)$, $n\in\N_0$, represent a control on the total number of individuals in the population at each generation. Thus, we distinguish three consecutive phases at each generation in this model: \emph{control phase}, \emph{mating phase} and \emph{reproduction phase}. The first phase is a control stage where the number of individuals that participate in the following phase is determined. This number is denoted as $\varphi_n(X_n)$ at generation $n$. Each one of these individuals can be female or male. It will be a female with probability $\lambda$ and a male with probability $1-\lambda$. The second phase is the mating phase where $\bar{F}_n$ females and $\bar{M}_n$ males at generation $n$ mate to produce $Y_n$ couples considering the promiscuous mating function as the mating function $L$, which is defined as $L(x,y)=x\min\{1,y\}$. The last phase is a reproduction phase where the couples give birth to their offspring. Thus, the variable $\xi_{ni}$ denotes the number of offspring of the $i$-th couple in the $n$-th generation. The probability distribution of this variable is named the \textit{offspring distribution} or \emph{reproduction law} and it is denoted $\{q_k\}_{k\in\N_0}$. To avoid trivialities it is assumed to satisfy $q_0+q_1+q_2< 1$.

\begin{remark}
\begin{enumerate}[label=(\alph*),ref=(\alph*)]
\item Condition \ref{cond:A2} is a natural hypothesis and it establishes that the number of couples is always less than or equal to the number of females. The meaning of the condition \ref{cond:A3} is obvious, if there is no female or male in the population, no couple can be formed.

\item Apart from the promiscuous mating other classic examples of mating functions used in two-sex branching processes theory and satisfying conditions \ref{cond:A2}-\ref{cond:A3} are the following:
\begin{enumerate}[label=(\roman*),ref=(\roman*)]
\item Polygamous mating: in this situation females and males could form couples with more than one individual of the opposite sex. It can be modelled through the function $L(x,y)=\min\{x,dy\}, \ d>1$.
\item \emph{Perfect fidelity mating}: in this case each individual mates with only one individual of the opposite sex forming an exclusive couple. The mating function is $L(x,y)=\min\{x,y\}$. 
\end{enumerate}
We note that although the results in this paper are provided for the promiscuous mating, the majority of them can be easily adapted for these other mating functions.
\end{enumerate}
\end{remark}

By the definition of the model, it is not difficult to check that the process $\{Y_n\}_{n\in\N_0}$ is a discrete homogeneous Markov process whose states are non-negative integers. Moreover, if in some generation there are not any couples, that is, $Y_n=0$ for certain $n>0$, and assuming that $\varphi_0(0)=0$ a.s. then, from that generation on, there will be neither individuals nor couples, i.e. $X_k=0$ and $Y_k=0$ for all $k>n$. This implies that the state 0 is absorbing and also the extinction of the population. Similar results about the classification of the states, the extinction and the asymptotic behaviour to those for the controlled branching process (see \cite[Chapters 3 and 4]{CBPs}) can be obtained. In the next results we only provide those that will be useful for the proofs of the results in the following sections. The proof of the first proposition is easily obtained with standard procedures and it is omitted. 

\begin{proposition}\label{thm:C2SBP-states-duality}
Let $\{Y_n\}_{n\in\N_0}$ be a BBPCI. Then:
\begin{enumerate}[label=(\roman*),ref=\emph{(\roman*)}]
\item $\mathbb{N}$ is a family of transient states.\label{thm:C2SBP-states-duality-i}
\item The classic duality extinction-explosion holds, that is, $P(Y_n\to 0)+P(Y_n\to \infty)=1$.\label{thm:C2SBP-states-duality-ii}
\end{enumerate}
\end{proposition}

In order to establish the following results let us denote the mean and variance of the variable $\xi_{01}$ as $0<\mu<\infty$, and $0<\delta^2<\infty$, respectively, and write $\epsilon(k)=E[\varphi_0(k)]$ and $\delta^2(k)=\V[\varphi_0(k)]$, and assume that both of them are positive and finite for each $k\in\N_0$.

\begin{theorem}\label{thm:C2SBP-suff-cond-extinc}
Let $\{Y_n\}_{n\in\N_0}$ be a BBPCI. If $\lambda\mu\epsilon(k)\leq k$, for each $k\in\N_0$, then $P(Y_n\to\infty|Y_0=y)=0$, for each $y\in\N$.
\end{theorem}

An immediate consequence of the previous theorem is that a controlled two-sex branching process with control functions $\varphi_0(k)$ following binomial distributions with parameters $k$ and $0<\rho<1$ becomes extinct a.s. if $\lambda\mu\rho\leq 1$. This remark will be useful in the proof of the results in Section~\ref{sec:Non-isolated}. 
 
\begin{theorem}\label{thm:C2SBP-nec-cond-extinc}
Let $\{Y_n\}_{n\in\N_0}$ be a BBPCI.
\begin{enumerate}[label=(\roman*),ref=\emph{(\roman*)}]
\item For each $n\in\N_0$, the conditional expectation is:\label{thm:C2SBP-nec-cond-extinc-i}
$$E[Y_{n+1}|Y_n]=\lambda E[\epsilon(X_{n+1})|Y_n]-\lambda E[h_{X_{n+1}}'(\lambda)|Y_n]\quad a.s.,$$
where $h_k(\cdot)$ denotes the probability generating function (p.g.f.) of $\varphi_0(k)$.

\item If the control function $\varphi_0(k)$ follows a binomial distribution with parameters $k$ and $\rho$, then:\label{thm:C2SBP-nec-cond-extinc-ii}
\begin{enumerate}[label=(\alph*),ref=\emph{(\alph*)}]
\item If $f(\cdot)$ denotes the p.g.f of the variable $\xi_{01}$, then there exist constants $C_1,C_2,C_3>0$ such that for $n\in\N_0$,\label{thm:C2SBP-nec-cond-extinc-ii-a}
\begin{align*}
E[Y_{n+1}|Y_n]&=\lambda\rho\mu Y_n - C_1  Y_n  f(1-\rho+\rho\lambda)^{ Y_n}\quad a.s.,\\
\V[Y_{n+1}|Y_n]&\leq C_2 Y_n  +C_3 Y_n^2 f(1-\rho+\rho\lambda)^{ Y_n}\quad a.s.
\end{align*}
\item If $\lambda\mu\rho> 1$, then $P(Y_n\to\infty|Y_0=y)>0$, for each $y\in \N$.\label{thm:C2SBP-nec-cond-extinc-ii-b}
\end{enumerate}
\end{enumerate}
\end{theorem}

\section{Definition of a predator-prey two-sex branching model}\label{sec:Definition}

Having described the controlled two-sex branching process, in this section we introduce a predator-prey two-sex branching process in order to model a predator-prey system.  We aim at introducing a model for a biological system where two animal species live together in the same environment, where one of them is the prey and the other one is its natural predator. We focus on the case that both species have sexual reproduction and  propose a controlled two-sex branching process to model the evolution of each of them, where the control mechanism is introduced in order to describe their natural interaction. We shall start with the formal definition of the model.

\vspace{2ex}

Let $\{t_{ni}: n\in\N_0, i\in\N\}$, $\{\tilde{t}_{ni}: n\in\N_0, i\in\N\}$, $\{\phi_{n}(t,\tilde{t}): n,t,\tilde{t}\in\N_0\}$ and $\{\tphi_{n}(t,\tilde{t}): n,t,\tilde{t}\in\N_0\}$ be independent families of non-negative  and integer valued r.v.s defined on the same probability space $(\Omega,\mathcal{A},P)$, and assume that:
\begin{enumerate}[label=\emph{(\roman*)},ref=\emph{(\roman*)}]
\item The r.v.s of the sequence $\{t_{ni}: n\in\N_0, i\in\N\}$ are i.i.d. with probability distribution $p=\{p_k\}_{k\in \mathbb{N}_0}$, where $p_k=P(t_{01}=k)$, $k\in\N_0$.

\item The r.v.s of the family $\{\tilde{t}_{ni}: n\in\N_0, i\in\N\}$ are i.i.d. with probability distribution $\tilde{p}=\{\tilde{p}_k\}_{k\in \mathbb{N}_0}$, where $\tilde{p}_k=P(\tilde{t}_{01}=k)$, $k\in\N_0$.

\item For each $n,t,\tilde{t}\in\N_0$, $\phi_{n}(t,\tilde{t})$ is a r.v. following a binomial distribution with parameters $t$ and $s(\tilde{t})$, where $s:\mathbb{R}\to[0,1]$ is a strictly increasing function which is continuous at 0 and such that, for certain $0\leq \rho_1< \rho_2\leq 1$,\label{cond:control-predator}
\begin{equation}\label{equ:lims}
s(0)=\rho_1, \quad \mbox{and} \quad \lim_{\tilde{t}\to \infty}s(\tilde{t})=\rho_2.
\end{equation}

\item For each $n,t,\tilde{t}\in\N_0$, $\tphi_{n}(t,\tilde{t})$ is a r.v. following a binomial distribution with parameters $\tilde{t}$ and $\ts(t)$, where $\ts:\mathbb{R}\to[0,1]$ is a strictly decreasing function which is continuous at 0 and such that, for certain $0\leq \trho_1< \trho_2\leq 1$, \label{cond:control-prey}
\begin{equation}\label{equ:limst}
\ts(0)=\trho_2, \quad \mbox{and} \quad \lim_{t\to \infty}\ts(t)=\trho_1.
\end{equation}
\end{enumerate}

A \textit{predator-prey two-sex branching process} is a bivariate stochastic process $\{(Z_n,\tZ_n)\}_{n\in \mathbb{N}_0}$ defined  recursively as:
\begin{equation}\label{def:model-total-indiv}
(Z_0,\tZ_0)=(N,\tilde{N}) \in \mathbb{N}^2,\qquad (T_{n+1},\tT_{n+1})=\left(\sum_{i=1}^{Z_{n}}t_{ni},\sum_{i=1}^{\tZ_{n}}\tilde{t}_{ni}\right),\qquad n\in\N_0,
\end{equation}
and 
\begin{equation*}
(Z_{n+1},\tZ_{n+1})=\left(L(F_{n+1},M_{n+1}),\tilde{L}(\tF_{n+1},\tM_{n+1})\right),\qquad n\in\N_0,
\end{equation*}
where conditionally on $\{\phi_{n+1}(T_{n+1},\tT_{n+1})=k\}$, the random vector $(F_{n+1},M_{n+1})$  follows a multinomial distribution with parameters $k$ and $(\alpha,1-\alpha)$, with $0<\alpha<1$, and conditionally on $\{\tphi_{n+1}(T_{n+1},\tT_{n+1})=\tilde{k}\}$ the random vector $(\tF_{n+1},\tM_{n+1})$ follows a multinomial distribution with parameters $\tilde{k}$ and $(\tilde{\alpha},1-\tilde{\alpha})$, with $0<\talpha<1$. Moreover, the empty sums in \eqref{def:model-total-indiv} are assumed to be 0.

\vspace{2ex}

The process defined above enables us to model a predator-prey system with non-overlapping generations. This process evolves as a three-stage procedure of reproduction, control and mating in each generation, where the previous variables have the next interpretation. The variables $T_n$ and $\tT_n$ represent the total number of predators and preys, respectively, at generation $n$, whereas $F_n$ and $M_n$ ($\tF_n$ and $\tM_n$) are the total numbers of progenitor predator females and predator males (prey females and prey males), respectively, at generation $n$. Moreover, $Z_n$ and $\tZ_n$ denote the number of predator couples and prey couples at the $n$-th generation, respectively. The dynamics of the three phases is described below.

In the \textit{reproduction phase}, couples of each species produce offspring independent of each others and in accordance with an offspring law. The offspring law may vary for the different species, but it remains constant over the generations for each species. Formally, the number of offspring produced by a couple of each species are represented by sequences of i.i.d. $\N_0$-valued r.v.s $\{t_{ni}: n\in\N_0, i\in\N\}$, and $\{\tilde{t}_{ni}: n\in\N_0, i\in\N\}$, where $t_{ni}$ denotes the number of offspring of the $i$-th predator couple at generation $n$ while $\tilde{t}_{ni}$ denotes the number of offspring of the $i$-th prey couple at generation $n$. The common probability distributions of these variables, $p=\{p_k\}_{k\in \mathbb{N}_0}$ and $\tilde{p}=\{\tilde{p}_k\}_{k\in \mathbb{N}_0}$, respectively, are called \textit{offspring distribution} or \emph{reproduction law of the predator and prey population}, respectively and we assume that they have finite and positive mean and variance, which are denoted $m$ and $\sigma^2$, respectively, for the predators and $\tilde{m}$ and $\tilde{\sigma}^2$ for the preys. We recall that we consider that the mating functions $L$ and $\tilde{L}$ are the promiscuous mating function and hence, in order to avoid trivialities, we also assume that $p_0+p_1+p_2<1$ and $\tilde{p}_0+\tilde{p}_1+\tilde{p}_2<1$. At the end of the reproduction phase at generation $n+1$, the sum of all the offspring of each species gives us the total number of predators, $T_{n+1}$, and the total number of preys, $\tT_{n+1}$, at this generation.

The reproduction stage is followed by the \textit{control phase}, where the number of predators and preys could be reduced due to several reasons such as their death because of the hunting, the lack of food supply or their capture by predators. Thus, if there are $T_{n+1}=t$ predators and $\tT_{n+1}=\tilde{t}$ preys in the population, the number of predators and preys which survive are given by the r.v.s $\phi_{n+1}(t,\tilde{t})$ and $\tilde{\phi}_{n+1}(t,\tilde{t})$, respectively. Considering that the survival of each predator (prey) is independent of the survival of the remaining predators (preys), and the probability of survival is the same for all the individuals in the same population, then it is natural to assume that the distributions of the variables $\phi_{n+1}(t,\tilde{t})$ and $\tilde{\phi}_{n+1}(t,\tilde{t})$ are binomial distributions. 

More specifically, $\phi_{n+1}(t,\tilde{t})$ is assumed to follow a binomial distribution with parameters $t$ and $s(\tilde{t})$, where the $s(\tilde{t})$ represents the survival probability of a predator given that there are $\tilde{t}$ preys in the population. The condition on the monotonicity of the function $s:\mathbb{R}\to[0,1]$ means that the smaller the number of preys is, the smaller the probability of survival of the predators. Regarding the conditions in \eqref{equ:lims}, we note that $\rho_1>0$ means that the predators could survive although there is no prey in the population because they could find another food source (other prey species, for instance). However, $\rho_1=0$ implies the extinction of the predator population when there are not any preys in the population at some generation. Moreover, $\rho_2<1$ means that, although there are enough preys in the population, the predators could die by several reasons (for example, by  hunting or their own predators), whereas $\rho_2=1$ indicates the survival of all the predators when the number of preys goes to infinity.  

Analogously, we assume that $\tphi_{n+1}(t,\tilde{t})$ follows a binomial distribution with parameters $\tilde{t}$ and $\ts(t)$, with $\ts(t)$ representing the survival probability of a prey given that there are $t$ predators in the population. The condition on the monotonicity of the function $\ts:\mathbb{R}\to[0,1]$ for the probability of survival of a prey indicates that the greater the number of predators is, the smaller the probability of survival of the preys becomes. The assumptions in \eqref{equ:limst} also have an intuitive interpretation. The condition $\trho_1>0$ means that the preys could survive although the number of predators goes to infinity because predators have limited appetite. The opposite situation, when $\trho_1=0$, implies the extinction of the prey population when there are a huge number of predators in the population at some generation. In addition, $\trho_2<1$ means that, despite the absence of predators in the ecosystem, the preys could die by another reasons (hunting or another predators). This case is excluded if $\trho_2=1$, which leads to the survival of all preys when there are not any predators in the environment in some generation.

At the end of this control phase, there are $F_{n+1}$ females and $M_{n+1}$ males within the survivor predator population at generation $n+1$. Thus, if $\alpha$ denotes the probability that a survivor predator is female, then the random vector $(F_{n+1},M_{n+1})$ follows a multinomial distribution of parameters $y$ and  $(\alpha,1-\alpha)$, given that $T_{n+1}=t$, $\tT_{n+1}=\tilde{t}$ and $\phi_{n+1}(t,\tilde{t})=y$. Similarly, there are $\tF_{n+1}$ females and $\tM_{n+1}$ males within the survivor prey population at generation $n+1$, and consequently, $(\tF_{n+1},\tM_{n+1})$ follows a multinomial distribution of parameters $\tilde{y}$, and $(\tilde{\alpha},1-\tilde{\alpha})$, given that $T_{n+1}=t$, $\tT_{n+1}=\tilde{t}$ and $\tilde{\phi}_{n+1}(t,\tilde{t})=\tilde{y}$, and where $\tilde{\alpha}$ is the probability that a survivor prey is female.

The last step is the \textit{mating phase}, where the predator and prey couples at generation $n+1$, $Z_{n+1}$ and $\tZ_{n+1}$, are determined by means of promiscuous mating functions depending on the number of females and males of each species at the current generation.

\begin{remark}\label{rem:functions}
One can propose several functions $s(\cdot)$ and $\tilde{s}(\cdot)$ satisfying the previous assumptions. For example, for the survival of the predator population
$$s(x)=\rho_2\left(1-a^{-x}\right)+\rho_1a^{-x},\quad a>1,\qquad s(x)=\rho_2\frac{x^k+\rho_1/\rho_2}{x^k+1},\quad k>0,$$
and similarly, for the survival of the prey population
$$\tilde{s}(x)=\left(\trho_2-\trho_1\right)a^{-x}+\trho_1,\quad a>1,\qquad \tilde{s}(x)=1-\left(1-\trho_1\right)\frac{x^k+\left(1-\trho_2\right)/\left(1-\trho_1\right)}{x^k+1},\quad k>0.$$
\end{remark}

\section{Basic properties of the model}\label{sec:BasicProperties}

In this section, we establish some basic properties of the process regarding the classification of the states and its main moments. First of all, note that from the definition of the model it is not difficult to deduce that the number of predators and prey couples in a certain generation only depend on the total number of predators and prey couples in the previous generation. Thus, the bivariate sequence $\{(Z_n,\tZ_n)\}_{n\in\N_0}$ is a discrete time homogeneous Markov chain whose states are two-dimensional vectors with non-negative integer coordinates. Moreover, it is immediate that $(0,0)$ is an absorbing state taking into account \ref{cond:control-predator} and \ref{cond:control-prey} and condition \ref{cond:A3}.

In the following easy-to-prove proposition we state some properties of the states associated with the process.

\begin{proposition}\label{prop:states}
Let $\{(Z_n,\tZ_n)\}_{n\in\N_0}$ be a predator-prey two-sex branching process. Then, 
\begin{enumerate}[label=(\roman*),ref=\emph{(\roman*)}]
\item Every non-null state $(i,j)\neq (0,0)$ is transient.\label{prop:states-i}
\item If $p_0+p_1+p_2<1$, $\tilde{p}_0+\tilde{p}_1+\tilde{p}_2<1$, and $0<\rho_1<\rho_2<1$, then the sets $\{(i,0): i>0\}$, $\{(0,j): j>0\}$ and $\{(i,j): i,j>0\}$ are classes of communicating states and each state leads to the state $(0,0)$. Furthermore, the process can move from the last set to the others in one step.\label{prop:states-promiscuous-iii}
\item If $p_0+p_1+p_2<1$, $\tilde{p}_0+\tilde{p}_1+\tilde{p}_2<1$ and $0=\rho_1<\rho_2=1$, then the sets $\{(0,j): j>0\}$ and $\{(i,j): i,j>0\}$ are classes of communicating states and each state leads to the state $(0,0)$. Furthermore, the states belonging to the second set may move to the other one and to the set $\{(i,0): i>0\}$ in one step. Finally,  the process moves from the last set to the state $(0,0)$ in one step.\label{prop:states-promiscuous-v}
\end{enumerate}
\end{proposition}

\vspace{2ex}

Next, we provide some results concerning the conditional moments of the variables involved in the definition of the process which will be useful in Section \ref{sec:Non-isolated}. Note that from the definition of the model, it is immediate to get the mean and variance of the control variables. Indeed, 
for $n,x,\tilde{x}\in\N_0$, the conditional expectations of the control variables are
\begin{align*}
\varepsilon(x,\tx)&:=E[\phi_{n}(T_n,\tT_n)|T_n=x,\tT_n=\tx]=xs(\tx),\\
\tilde{\varepsilon}(x,\tx)&:=E[\tphi_{n}(T_n,\tT_n)|T_n=x,\tT_n=\tx]=\tx\ts(x).
\end{align*}
and the conditional variances of the control variables are
\begin{align*}
\sigma^2(x,\tx)&:=Var[\phi_{n}(T_n,\tT_n)|T_n=x,\tT_n=\tx]=xs(\tx)(1-s(\tx)),\\
\tilde{\sigma}^2(x,\tx)&:=Var[\tphi_{n}(T_n,\tT_n)|T_n=x,\tT_n=\tx]=\tx\ts(x)(1-\ts(x)).
\end{align*}



For the next results, we introduce the following notation concerning the $\sigma$-algebras generated by the variables involved in the definition of the model. In particular, we denote
\begin{align*}
\mathcal{F}_n&=\sigma(Z_l,\tZ_l:l=0,\ldots,n),\quad n\in\N_0,\\
\mathcal{G}_n&=\sigma(Z_l,\tZ_l,T_{l+1},\tilde{T}_{l+1}:l=0,\ldots,n-1),\quad n\in\N,\\
\mathcal{H}_{n}&=\sigma(Z_l,\tZ_l,T_{l+1},\tilde{T}_{l+1},\phi_{l+1}(T_{l+1},\tilde{T}_{l+1}),\tphi_{l+1}(T_{l+1},\tilde{T}_{l+1}):l=0,\ldots,n-1),\quad n\in\N.
\end{align*}
Then, we have that $\mathcal{F}_{n-1}\subseteq \mathcal{G}_{n}\subseteq\mathcal{H}_{n}$, for $n\in\N$. Now, the conditional moments of the number of predator and prey individuals can be easily obtained from the definition of the model and hence, the proof of the next proposition is omitted.

\begin{proposition}[Conditional moments of the number of individuals]\label{Moments Individuals} 
Let $\{(Z_n,\tZ_n)\}_{n\in\N_0}$ be a predator-prey two-sex branching process. Then, for each $n\in\N_0$,
\begin{enumerate}[label=(\roman*),ref=\emph{(\roman*)}]
\item $E[T_{n+1}|\mathcal{F}_n]=m Z_n$, \quad and \quad $E[\tT_{n+1}|\mathcal{F}_n]=m\tZ_n$.\label{Moments Individuals-i} 
\item $Var[T_{n+1}|\mathcal{F}_n]=\sigma^2 Z_n$, \quad and \quad $Var[\tT_{n+1}|\mathcal{F}_n]=\tilde{\sigma}^2\tZ_n$.\label{Moments Individuals-ii} 
\end{enumerate}
\end{proposition}

Next, we establish some results concerning to the conditional moments of the total number of female and male predators and female and male preys.

\begin{proposition}[Conditional moments of the number of females and males]\label{Moments FM|T}
Let \linebreak$\{(Z_n,\tZ_n)\}_{n\in\N_0}$ be a predator-prey two-sex branching process. Then, for each $n\in\N_0$,
\begin{enumerate}[label=(\roman*),ref=\emph{(\roman*)}]
\item The expected number of predator and prey females given the number of individuals and couples are\label{Moments FM|T-i}
$$E[F_{n}|\mathcal{G}_n]=\alpha T_n s(\tT_n), \quad \text{and} \quad E[\tF_{n}|\mathcal{G}_n]=\talpha \tT_n\ts(T_n),$$
and the corresponding to the number of predator and prey males are
$$E[M_{n}|\mathcal{G}_n]=(1-\alpha)T_n s(\tT_n), \quad \text{and}  \quad E[\tM_{n}|\mathcal{G}_n]=(1-\talpha)\tT_n\ts(T_n).$$
\item The conditional variances of the number of predator and prey females and males given the number of individuals and couples are \label{Moments FM|T-ii}
\begin{align*}
Var[F_{n}|\mathcal{G}_n]&=\alpha^2 T_n s(\tT_n)(1-s(\tT_n))+\alpha(1-\alpha)T_n s(\tT_n),\\
Var[\tF_{n}|\mathcal{G}_n]&=\talpha^2 \tT_n\ts(T_n)(1-\ts(T_n))+\talpha(1-\talpha)\tT_n\ts(T_n),
\end{align*}
and
\begin{align*}
Var[M_{n}|\mathcal{G}_n]&=(1-\alpha)^2 T_n s(\tT_n)(1-s(\tT_n))+\alpha(1-\alpha)T_n s(\tT_n),\\ Var[\tM_{n}|\mathcal{G}_n]&=(1-\talpha)^2 \tT_n\ts(T_n)(1-\ts(T_n))+\talpha(1-\talpha)\tT_n\ts(T_n),
\end{align*}
respectively.
\end{enumerate}
\end{proposition}

\begin{proposition}[Conditional moments of the number of females and males]\label{Moments FM|Z}
Let\linebreak $\{(Z_n,\tZ_n)\}_{n\in\N_0}$ be a predator-prey two-sex branching process. Then, for each $n\in\N_0$,
\begin{enumerate}[label=(\roman*),ref=\emph{(\roman*)}]
\item The expected number of predator and prey females given the number of couples are\label{Moments FM|Z-i}
$$E[F_{n+1}|\mathcal{F}_n]=\alpha m Z_{n}E[s(\tT_{n+1})|\mathcal{F}_n],\quad\text{ and } \quad E[\tF_{n+1}|\mathcal{F}_n]=\tilde{\alpha} \tilde{m} \tZ_{n}E[\ts(T_{n+1})|\mathcal{F}_n],$$ 
and the corresponding to the number of predator and prey males are
$$E[M_{n+1}|\mathcal{F}_n]=(1-\alpha) m Z_{n}E[s(\tT_{n+1})|\mathcal{F}_n],\quad\text{ and }\quad E[\tM_{n+1}|\mathcal{F}_n]=(1-\tilde{\alpha}) \tilde{m} \tZ_{n}E[\ts(T_{n+1})|\mathcal{F}_n].$$

\item The conditional variance of the number of predator and prey females and males given the number of couples are\label{Moments FM|Z-ii}
{\footnotesize\begin{align*}
Var[F_{n+1}|\mathcal{F}_n]&=\alpha^2 m^2 Z_n^2Var[s(\tT_{n+1})|\mathcal{F}_n]+\alpha m Z_n E[s(\tT_{n+1})|\mathcal{F}_n]+ \alpha^2 Z_n (\sigma^2-m)E[s(\tT_{n+1})^2|\mathcal{F}_n],\\
Var[\tF_{n+1}|\mathcal{F}_n]&=\talpha^2 \tilde{m}^2 \tZ_n^2Var[\ts(T_{n+1})|\mathcal{F}_n]+\talpha \tilde{m} \tZ_n E[\ts(T_{n+1})|\mathcal{F}_n]+ \talpha^2 \tZ_n (\tilde{\sigma}^2-\tilde{m})E[\ts(T_{n+1})^2|\mathcal{F}_n].
\end{align*}}
and 
\begin{align*}
Var[M_{n+1}|\mathcal{F}_n]&=(1-\alpha)^2 m^2 Z_n^2Var[s(\tT_{n+1})|\mathcal{F}_n]+(1-\alpha) m Z_n E[s(\tT_{n+1})|\mathcal{F}_n]\\
&\phantom{=}+ (1-\alpha)^2 Z_n (\sigma^2-m)E[s(\tT_{n+1})^2|\mathcal{F}_n],\\
Var[\tM_{n+1}|\mathcal{F}_n]&=(1-\talpha)^2 \tilde{m}^2 \tZ_n^2Var[\ts(T_{n+1})|\mathcal{F}_n]+(1-\talpha) \tilde{m} \tZ_n E[\ts(T_{n+1})|\mathcal{F}_n]\\
&\phantom{=}+ (1-\talpha)^2 \tZ_n (\tilde{\sigma}^2-\tilde{m})E[\ts(T_{n+1})^2|\mathcal{F}_n],
\end{align*}
respectively.
\end{enumerate}
\end{proposition}

We note that from this result, it is easy to obtain an upper bound for the conditional expectation of the number of predator and prey couples as stated in the following proposition. 
 
\begin{proposition}\label{prop:expectation}
Let $\{(Z_n,\tZ_n)\}_{n \in \mathbb{N}_0}$ be a predator-prey two-sex branching process. Then, for each $n\in\N_0$,
$$E[Z_{n+1}|\mathcal{F}_n]\leq\alpha m Z_n E[s(\tT_n)|\mathcal{F}_n]\leq
\rho_2\alpha m Z_n, \quad \mbox{and} \quad E[\tZ_{n+1}|\mathcal{F}_n]\leq\tilde{\alpha}\tilde{m}\tZ_n E[\tilde{s}(T_n)|\mathcal{F}_n]\leq
\trho_2\tilde{\alpha}\tilde{m}\tZ_n.$$
\end{proposition}

\begin{remark}
Since $\rho_1\leq s(t)\leq \rho_2$, and $\trho_1\leq \ts(t)\leq \trho_2$, for all $t\geq 0$, we obtain upper and lower bounds for the previous conditional expectations,
\begin{eqnarray*}
\rho_1 \alpha T_{n}\leq &E[F_{n}|\mathcal{G}_n]&\leq \rho_2 \alpha T_{n},\\
\rho_1 \alpha m Z_{n}\leq &E[F_{n+1}|\mathcal{F}_n]&\leq \rho_2 \alpha m Z_{n}.
\end{eqnarray*}
Similar arguments let us get the bounds for the remaining conditional expectations.
\end{remark}

Now, we provide the conditional moments of the number of predator and prey individuals given the total number of individuals and couples in previous generations.

\begin{proposition}[Conditional moments of the number of individuals]\label{lem:moments-T-promiscuous}
Let $\{(Z_n,\tZ_n)\}_{n\in\N_0}$ be a predator-prey two-sex branching process. Then, for any $n\in\N$,
\begin{enumerate}[label=(\roman*),ref=\emph{(\roman*)}]
\item $E[T_{n+1}|\mathcal{G}_n]=\alpha m T_ns(\tT_n)-\alpha m T_ns(\tT_n)(1-s(\tT_n)+\alpha s(\tT_n))^{T_n-1}$.\\
 $E[\tT_{n+1}|\mathcal{G}_n]=\talpha \tm \tT_n\ts(T_n)-\talpha \tm \tT_n\ts(T_n)(1-\ts(T_n)+\talpha \ts(T_n))^{\tT_n-1}$.\label{lem:moments-T-promiscuous-i}

\item $Var[T_{n+1}|\mathcal{G}_n]\leq (\sigma^2+m^2)\alpha T_ns(\tT_n)+2\alpha^2 m^2 T_n^2 s(\tT_n)^2 (1-s(\tT_n)+\alpha s(\tT_n))^{T_n-1}.$\\
$Var[\tT_{n+1}|\mathcal{G}_n]\leq (\tilde{\sigma}^2+\tm^2)\talpha\tT_n \ts(T_n)+2\talpha^2 \tm^2 \tT_n^2 \ts(T_n)^2 (1-\ts(T_n)+\talpha \ts(T_n))^{\tT_n-1}.$\label{lem:moments-T-promiscuous-ii}
\end{enumerate}
\end{proposition}

\vspace{2ex}

To conclude this section we establish the following result where we derive the usual property of branching processes known as the extinction-explosion dichotomy. 

\begin{proposition}\label{prop:extexp} Let $\{(Z_n,\tZ_n)\}_{n\in\N_0}$ be a predator-prey two-sex branching process. Then,
\begin{enumerate}[label=(\roman*),ref=\emph{(\roman*)}]
\item $P(\liminf_{n\to\infty}(Z_n,\tZ_n)=(k,\tilde{k}))=0$, and $P(\limsup_{n\to\infty}(Z_n,\tZ_n)=(k,\tilde{k}))=0$, for each $(k,\tilde{k})\in\N_0^2\backslash\{(0,0)\}$.\label{prop:extexp-i}
\item $P(Z_n\to 0,\tZ_n\to 0)+ P(Z_n \to\infty,\tZ_n\to \infty)+P(Z_n\to \infty,\tZ_n\to 0)+
P(Z_n \to 0,\tZ_n\to\infty)=1$.\label{prop:extexp-ii}
\end{enumerate}
\end{proposition}

The sets $\{Z_n\to 0,\tZ_n\to 0\}$, $\{Z_n \to\infty,\tZ_n\to \infty\}$, $\{Z_n\to \infty,\tZ_n\to 0\}$ and $\{Z_n \to 0,\tZ_n\to\infty\}$ are termed extinction of both populations, survival (or coexistence) of both populations, predator population fixation, and prey population fixation, respectively. Moreover, if we denote the extinction and survival of the predator population as $\{Z_n\to 0\}$, and $\{Z_n\to \infty\}$, and the extinction and survival of the prey population as $\{\tZ_n\to 0\}$, and $\{\tZ_n\to \infty\}$, then in view of Proposition \ref{prop:extexp} it is immediate that,
\begin{align}
\{Z_n\to 0\}&=\{Z_n\to 0, \tZ_n\to 0\} \cup \{Z_n\to 0, \tZ_n\to \infty\}\quad a.s.\nonumber\\
\{Z_n\to \infty\}&=\{Z_n\to \infty, \tZ_n\to 0\} \cup \{Z_n\to \infty, \tZ_n\to \infty\}\quad a.s.\label{equ:Zinfty}\\
\{\tZ_n\to 0\}&=\{Z_n\to 0, \tZ_n\to 0\} \cup \{Z_n\to \infty, \tZ_n\to 0\}\quad a.s.\nonumber\\
\{\tZ_n\to \infty\}&=\{Z_n\to 0, \tZ_n\to \infty\} \cup \{Z_n\to \infty, \tZ_n\to \infty\}\quad a.s.\label{equ:Ztinfty}
\end{align}

\begin{remark}
Note that, in the cases of fixation of the predator and the prey, the process behaves as a BBPCI defined \eqref{def:C2SBP} in Section \ref{sec:C2SBP}.
\end{remark}

\section{Predator-prey isolated system}\label{sec:Isolated}

In this section we consider an isolated predator-prey system, that is, we assume that both species live together in an isolated area where the prey population constitutes the only food resource for the predators. More specifically, we focus on the case that there is an autochthonous species living in an isolated region and a new (invasive) species is introduced in that ecosystem and it preys on the autochthonous one. The question of how those populations evolve together is tackled in this section and it is of great interest for the preservation of the species in these environments. Several examples of those geographically isolated regions have been reported such as, for instance, Azores Islands (see \cite{Moreira2019}), Eastern Island (see \cite{Sottornoff2019}), Macquarie Island (see \cite{Frugone2019}), isolated regions in Finland and Northwest Russia (see \cite{Vila2003}) or seafloor plateau (see \cite{Richter2019}). 

In terms of our model, this situation can be expressed as $\rho_1=0$, since the fact that the prey is the only food supply for the predators implies that the probability of survival of any predator is zero if there is no prey in the population. Analogously, we consider $\trho_1=0$ which means an unlimited appetite of the predators that implies a null probability of survival for all preys when the number of predators in the population goes to infinity. We also allow $\rho_2\leq 1$ and $\trho_2\leq 1$ because both the predators and preys could die due to natural causes although there is no prey or predator, respectively, in the ecosystem. 

Notice that, under these assumptions if the number of prey couples at some generation is equal to zero then the prey population becomes extinct forever, that is, if $\tZ_n=0$ for some $n>0$, then $\tT_k=0$ and $\tZ_k=0$ for all $k>n$. The extinction of the prey population also bounds the predator population to disappear due to the fact that $s(0)=0$, and consequently, $\phi_k(t,0)=0$, for all $k>n$, from which one deduces $Z_k=0$, and $T_{k+1}=0$, for all $k>n$. On the other hand, if at some generation $n$ there are no predator couples, then the predator population becomes extinct forever, i.e., $T_k=0$, and $Z_k=0$, for all $k>n$. As a result, and since $\ts(0)=\rho_2$, $\tilde{\phi}_k(\tilde{t},0)$ follows a binomial distribution with parameters $\tilde{t}$ and $\rho_2$, for all $k>n$, which means that the prey population behaves as a BBPCI (see Section \ref{sec:C2SBP}). 

Bearing in mind these considerations, we study the fate of both species in the population. To that end, given $i,j>0$, we write $P_{(i,j)}(\cdot)=P(\cdot|(Z_0,\tZ_0)=(i,j))$, in the remaining of this paper.

\subsection{The certain extinction of the predator population}\label{subsec:predator_fixation}

The following result is very natural from the definition of the model and it means that in this kind of populations we cannot have the extinction of the prey population and the survival of the predator population (predator fixation). Recall that the prey population is the only food supply of the predator population, and hence, the extinction of preys dooms the predator population to the extinction.

\begin{proposition}\label{prop:predator fixation}
Let $\{(Z_n,\tZ_n)\}_{n\in\N_0}$ be a predator-prey two-sex branching process. Then, for any initial values  $i,j>0$, $P_{(i,j)}(Z_n\to \infty, \tZ_n\to 0)=0$.
\end{proposition}


The following result establishes that in this kind of systems the coexistence of preys and predators is not possible. Intuitively, this is deduced as follows. When the number of predators is too large, the probability of survival of the preys is too low and consequently, we expect a huge drop in the number of preys and as a result, in the number of predators in the next generation. This cycle can be repeated until a generation where the probability of survival of each prey is so negligible that the entire population of preys becomes extinct, so does the predator population in the following generation.

\begin{proposition}\label{prop:impossible coexistence}
Let $\{(Z_n,\tZ_n)\}_{n\in\N_0}$ be a predator-prey two-sex branching process. Then, for any initial values $i,j>0$, $P_{(i,j)}(Z_n\to \infty, \tZ_n\to \infty)=0$.
\end{proposition}

Note that by \eqref{equ:Zinfty}, the previous results and Proposition \ref{prop:extexp}~\ref{prop:extexp-ii} we obtain the certain extinction of the predator population, that is, $P_{(i,j)}(Z_n\to 0)=1$, for any initial values $i,j>0$.

\subsection{The fixation of the prey population}

From Propositions \ref{prop:predator fixation} and \ref{prop:impossible coexistence}, we deduce that the predator population becomes extinct almost surely in this model. The question that arises is whether the prey population has a chance to survive depending on its reproductive capacity once the predator population has become extinct. Let us start with the case $\trho_2=1$. On the prey fixation set, the prey population behaves as a standard two-sex branching process (without any kind of control) from one generation on once the predator population has become extinct. Thus, the theory developed in \cite{daleyA} can be applied in this setting and we immediately obtain the following result.

\begin{proposition}\label{prop:prey-fixation-Daley}
Let $\{(Z_n,\tZ_n)\}_{n\in\N_0}$ be a predator-prey two-sex branching process. If $\trho_2=1$, for any initial values $i,j>0$, $P(Z_n\to 0,\tZ_n\to\infty)>0$ if and only if $\talpha \tm>1$.
\end{proposition}

On the other hand, if $\trho_2<1$, on the prey fixation set the prey population behaves as the BBPCI introduced in Section \ref{sec:C2SBP} from one generation on. The control variables in this model follow binomial distributions with size equal to the number of preys at the corresponding generations and probability $\trho_2$. In this case,  the result is a direct consequence of Theorems \ref{thm:C2SBP-suff-cond-extinc} and \ref{thm:C2SBP-nec-cond-extinc}~\ref{thm:C2SBP-nec-cond-extinc-ii}~\ref{thm:C2SBP-nec-cond-extinc-ii-b} provided in the mentioned section.

\begin{proposition}\label{prop:prey-fixation-promiscuous}
Let $\{(Z_n,\tZ_n)\}_{n\in\N_0}$ be a predator-prey two-sex branching process. If $\trho_2<1$, for any initial values $i,j>0$, $P_{(i,j)}(Z_n\to 0,\tZ_n\to\infty)>0$ if and only if $\talpha \tm\trho_2>1$.
\end{proposition}

\section{Predator-prey non-isolated system}\label{sec:Non-isolated}

In this section we consider the case that $0<\rho_1<\rho_2\leq 1$ and $0<\trho_1<\trho_2<1$. Thus, contrary to Section \ref{sec:Isolated}, the predators have a positive probability of survival in absence of the prey population ($\rho_1>0$) due to the availability of other food resources and also a limited appetite ($\trho_1>0$) which allows the prey population to have a positive probability of survival even when the predator population size goes to infinity. Moreover, individuals of the prey population might not reproduce because of the presence of predators, but also for other reasons such as their hunting, diseases, or migratory movements  ($\trho_2<1$).

\subsection{The fixation of the predator and prey populations}\label{sec:Non-isolated fixation}

In this subsection we study necessary and sufficient conditions for the fixation of each species, that is, for one of the two species (predator or prey) to survive and the other one to become extinct. In the fixation events, the surviving species behaves as the BBPCI introduced in Section \ref{sec:C2SBP} from some generation on. The corresponding offspring distribution is the reproduction law of the survivor species and the control functions follows binomial distributions with constant probability of success $\gamma$, where $\gamma=\rho_1$ in the case of the predator fixation and $\gamma=\trho_2$ in the case of the prey fixation. Thus, by using Theorems \ref{thm:C2SBP-suff-cond-extinc} and \ref{thm:C2SBP-nec-cond-extinc} \ref{thm:C2SBP-nec-cond-extinc-ii} \ref{thm:C2SBP-nec-cond-extinc-ii-b} we have the following result.

\begin{proposition}\label{prop: non isolated prey and predator fixation}
Let $\{(Z_n,\tZ_n)\}_{n\in\N_0}$ be a predator-prey two-sex branching process. For any initial values $i,j>0$:
\begin{enumerate}[label=(\roman*),ref=\emph{(\roman*)}]
\item $P_{(i,j)}(Z_n\to \infty,\tZ_n\to 0)>0$ if and only if $\rho_1\alpha m>1$.\label{prop:non-iso-predator-fixation}
\item $P_{(i,j)}(Z_n\to 0,\tZ_n\to \infty)>0$ if and only if $\trho_2\talpha \tilde{m}>1$.\label{prop:non-iso-prey-fixation}
\end{enumerate}
\end{proposition}

Intuitively, this result states that a necessary and sufficient condition for the predator population to have a positive probability of fixation is that the mean number of female predators which survive after the control is greater than one. Alternatively, the second part of the result also states that a necessary and sufficient condition for the prey population to have a positive probability of fixation is that the mean number of female preys which survive after the control is greater than one. Now, taking into account \eqref{equ:Zinfty} and \eqref{equ:Ztinfty}, an immediate consequence of Proposition \ref{prop: non isolated prey and predator fixation} is the following corollary:

\begin{corolario}\label{coro: non isolated prey and predator fixation}
Let $\{(Z_n,\tZ_n)\}_{n\in\N_0}$ be a predator-prey two-sex branching process. For any initial values $i,j>0$:
\begin{enumerate}[label=(\roman*),ref=\emph{(\roman*)}]
\item If $\rho_1\alpha m>1$, then $P_{(i,j)}(Z_n\to \infty)>0$.\label{coro: non isolated prey and predator fixation-i}
\item If $\trho_2\talpha \tilde{m}>1$, then $P_{(i,j)}(\tZ_n\to \infty)>0$.\label{coro: non isolated prey and predator fixation-ii}
\end{enumerate}
\end{corolario}

\subsection{The extinction of the population}

Intuitively, it is clear that if the mean number of female predators (preys) which survives after the control is less than one then the reproductive capacity of the species is not enough to keep it alive by its own. Therefore, we can establish the following result:

\begin{proposition}\label{prop:non-isol extinction}
Let $\{(Z_n,\tZ_n)\}_{n\in\N_0}$ be a predator-prey two-sex branching process. For any initial values $i,j>0$:
\begin{enumerate}[label=(\roman*),ref=\emph{(\roman*)}]
\item If $\rho_2\alpha m\leq 1$, then $P_{(i,j)}(Z_n\to 0)=1$.\label{prop:non-isol extinction-i}
\item If $\trho_2\talpha \tilde{m}\leq 1$, then $P_{(i,j)}(\tZ_n\to 0)=1$.\label{prop:non-isol extinction-ii}
\end{enumerate}
\end{proposition}

Taking into account \eqref{equ:Zinfty} and \eqref{equ:Ztinfty} and from Proposition  \ref{prop:non-isol extinction}, we deduce the following result on the coexistence of the species:

\begin{corolario}\label{coro:non-isol extinction}
Let $\{(Z_n,\tZ_n)\}_{n\in\N_0}$ be a predator-prey two-sex branching process. For any initial values $i,j>0$, if $\min\{\rho_2\alpha m, \trho_2\talpha \tilde{m}\}\leq 1$, then $P_{(i,j)}(Z_n\to \infty, \tZ_n\to \infty)=0$.
\end{corolario}

We note that there is always a positive probability for the complete extinction of the predator-prey system. This could happens for several reasons: either because  all individuals of both species might die during the control phase or because all the survivors of both populations might be of the same sex, which makes impossible to form new couples. If we also allow $p_0>0$ and $\tp_0>0$, then there is a positive probability that the predator and prey couples produce no offspring.

In the next result, we determine a necessary and sufficient condition for both species to become extinct with probability one, which means the extinction of the entire predator-prey system.

\begin{proposition}\label{prop:non-isol common extinction}
Let $\{(Z_n,\tZ_n)\}_{n\in\N_0}$ be a predator-prey two-sex branching process. For any initial values $i,j>0$, $P_{(i,j)}(Z_n\to 0, \tZ_n\to 0)=1$ if and only if $\max\{\rho_1m\alpha,\trho_2\tilde{m}\tilde{\alpha}\}\leq 1$. 
\end{proposition}

\subsection{The predator and prey coexistence}

In the following result, we study the possibility of having the coexistence of the predator and prey populations. 

\begin{theorem}\label{prop:non-isol coexistence}
Let $\{(Z_n,\tZ_n)\}_{n\in\N_0}$ be a predator-prey two-sex branching process. For any initial values $i,j>0$:
\begin{enumerate}[label=(\roman*),ref=\emph{(\roman*)}]
\item If $\trho_1\tilde{m}\tilde{\alpha}<1$, then $P_{(i,j)}(Z_n\to \infty,\tZ_n\to \infty)=0$.\label{prop:non-isol coexistence-i}
\item If $\min\{\rho_2 m\alpha,\trho_1\tilde{m}\tilde{\alpha}\}>1$, then $P_{(i,j)}(Z_n\to \infty, \tZ_n\to \infty)>0$. \label{prop:non-isol coexistence-ii}
\end{enumerate}
\end{theorem}


\vspace{0.25cm}

\begin{figure}
\centering
  \epsfig{file=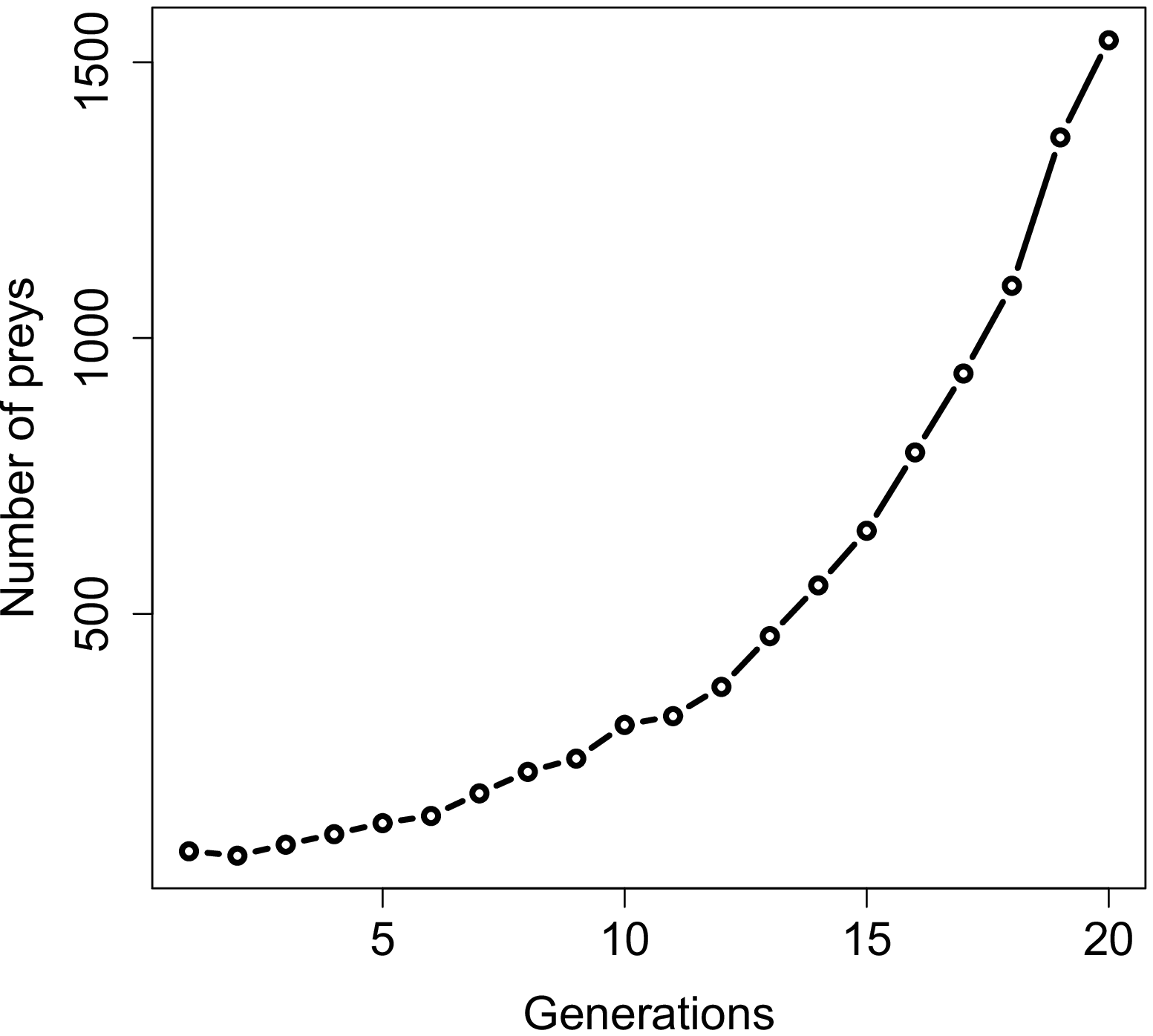, width=0.45\linewidth}\hspace{0.05\linewidth}
  \epsfig{file=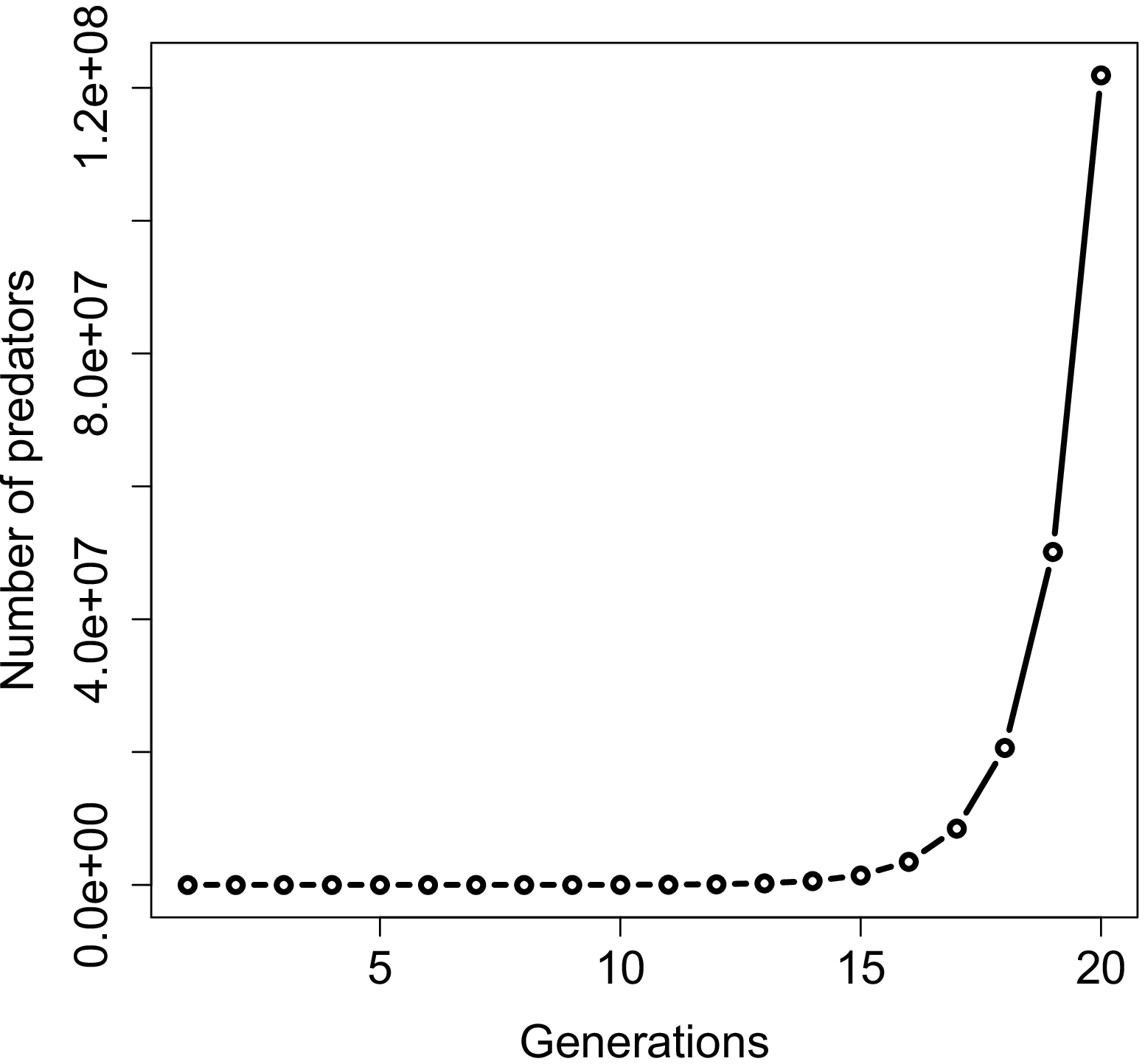, width=0.45\linewidth}
  \captionof{figure}{Evolution of the number of preys (left) and predators (right) during 20 generations of a predator-prey two sex branching process with parameters satisfying $\alpha m \rho_1> 1$ and $\talpha \tilde{m} \trho_1> 1$. The process started with $Z_0 = 1$ predator and $\tZ_0=50$ preys and the parameters of the model are $\alpha= 0.44$, $\talpha= 0.48$, $m = 6.5$, $\tilde{m}=7$, $\rho_1= 0.35$, $\rho_2 = 0.85$, $\trho_1=0.35$ and $\trho_2=0.98$. }\label{exponential}
  \end{figure}

In the case of coexistence, both populations grow geometrically at certain rate determined by the parameters of the model (see Figure \ref{exponential}). Note that this kind of growth is the classical behaviour in branching processes including the Galton-Watson predator-prey process (see \cite{Alsmeyer-1993}), but it is not typical in predator-prey systems modelled through ODEs, where periodic cycles are observed. However, there are examples of populations with this exponential growth behaviour (see \cite{Vila2003}, \cite{Stewart2011}, \cite{Omeja2014}, \cite{Wang2014} or \cite{Chapman2018}). Moreover, our model can be also applied at initial stages of other populations when there is a small number of individuals, such as, for example, in populations of endangered species where the exponential growth is shown (see \cite{Law2015}, \cite{Guevara2016} or \cite{Li2017}).

\section{Discussion}\label{sec:Discussion}

In this paper a controlled two-sex branching process is introduced with the aim of modelling predator-prey interactions in different situations. The generation-by-generation evolution of each species is also studied. We have focussed on the case where both populations have sexual reproduction and females and males form couples via a promiscuous sexual mating system. The process evolves in three stages which are repeated at each generation. First, couples of each species produce a random number of offspring (reproduction phase). Next, some of these individuals die because of the interaction between preys and predators and the remaining ones survive and are able to mate (control phase). This control mechanism is modelled through binomial distributions where for each species the probability of success depends on the number of individuals of the other one. Similarly, once that the number of survivors of each species is known, multinomial distributions are applied to determine the number of females and males among them, where the probability vector is the pair given by the probability that one individual is female and the probability that it is male. The final step at any generation is the mating phase when the couples that are able to produce offspring in the next generation are formed.

We highlight that the main novelty of our model with respect to previous ones is the sexual reproduction among the individuals within each species. This important feature had not been considered in the literature yet, but it is quite common in the nature. We also remark that although we have considered a promiscuous mating, other mating systems could be applied  and our results could be adapted to those cases modelling different situations.

To analyse the temporal evolution of the process, we have considered two possible scenarios. In the first one we have assumed that the preys live in an isolated region and a new species, the predator, is introduced. In the second one, both species live together in certain area where predators have different food sources (not only the prey) and a limit appetite; this enables the survival of preys even when there is a large number of predators. On the one hand, our results show that the prey population has two possible behaviours (extinction or explosion) in both settings and this behaviour depends on the mean number of prey females that survives during the control phase. On the other hand, the fate of the predator population depends on the situation. In the first scenario, the predator group becomes extinct almost surely because the prey is its only food supply and the predator has limitless appetite. In the second one, the predator could become extinct or growth indefinitely depending on the mean number of female predators that survive during the natural control of the population.   

As mentioned above, the model presented in this paper captures the fact that an increment (drop) in the number of preys implies an increment (drop) in the predator population size, but it does not describe the fluctuation behaviour where a large number of preys and predators is followed after some time by a smaller number of preys and predators (see the classical example of the snowshoe hare and the Canadian lynx in \cite{Gilpin1973}). These type of oscillations are usually caused by fluctuations in the environmental conditions or periodic changes in the reproductive capacity of the species. For instance, in the case of the snowshoe hare, the litter sizes vary between years and females only give birth during the breeding season, which is stimulated by new plants, and it begins around mid-March and runs to August (see \cite{Dodds-1965} or \cite{Krebs-Boonstra-2001}). Then, an appropriate way to introduce this periodic reproductive behaviour in the model is to let the predator offspring distribution and the prey offspring distribution change over the time.

Moreover, even in those populations where an exponential growth is observed, sooner or later the growth stops due to environmental conditions or other causes (for instance, the growth of a population of turtles truncated due to an oil spill is reported in \cite{Kocmoud2019}). If we wish to reflect the saturation of the environment -and cut the initial exponential growth- it is necessary to modify our assumptions.  For example, we could introduce a carrying capacity parameter in the probability functions $s(\cdot)$ and $\ts(\cdot)$ or even let them depend on two values simultaneously: the number of the preys and the number of the predators. The introduction of all the aforementioned modifications would lead to a more complex model which is beyond of the scope of this paper and it is left for future research.   

%
%

\section*{Acknowledgements}

This research was supported by the Ministerio de Econom{\'i}a y Competitividad and  the Fondo Europeo de Desarrollo Regional (MINECO/FEDER, UE) [grant MTM2015-70522-P] and by the Junta de Extremadura/European Regional Development Fund [grants IB16103 and GR18103].

The authors would like to thank Dr. Miguel Gonz\'alez and Dr. In\'es del Puerto (University of Extremadura) for their valuable discussions.

\appendix

\section*{Appendices}

\section{Proofs of the results in Section \ref{sec:C2SBP}}

In order to facilitate the reading of the proofs, in this section we denote $\bar{\mathcal{F}}_{n}=\sigma(Y_l:l=0,\ldots n)$, $\bar{\mathcal{G}}_{n}=\sigma(Y_l,X_{l+1}:l=0,\ldots,n-1)$, and $\bar{\mathcal{H}}_{n}=\sigma(Y_l,X_{l+1},\varphi_{l+1}(X_{l+1}):l=0,\ldots,n-1)$, and note that $\bar{\mathcal{F}}_{n-1}\subseteq \bar{\mathcal{G}}_n \subseteq \bar{\mathcal{H}}_n$, for $n\in\N$.

\vspace{0.5cm}

\begin{Prf}[Theorem \ref{thm:C2SBP-suff-cond-extinc}]
The proof follows the arguments of Theorem 1 in \cite{art-2002}. Since
\begin{align*}
E[Y_{n+1}|\bar{\mathcal{F}}_{n}] 
&\leq E[E[\bar{F}_{n+1}|\bar{\mathcal{H}}_{n+1}]|\bar{\mathcal{F}}_{n}]\\
&=\lambda E\left[E[\varphi_{n+1}(X_{n+1})|\bar{\mathcal{G}}_{n+1}]|\bar{\mathcal{F}}_{n}\right]\\
&\leq E\left[X_{n+1}|\bar{\mathcal{F}}_{n}\right]\mu^{-1}
=Y_n\quad a.s.,
\end{align*}
the sequence $\{E[Y_{n}]\}_{n\in\N_0}$ is decreasing and bounded from above by $E[Y_0]$. Thus, by Fatou's lemma
$$E[\liminf_{n\to\infty}Y_n]\leq \liminf_{n\to\infty} E[Y_n]<\infty,$$
and hence, $\liminf_{n\to\infty}Y_n$ is finite a.s., from which we get $P(Y_n\to\infty|Y_0=y)=0$.
\end{Prf}

\begin{Prf}[Theorem \ref{thm:C2SBP-nec-cond-extinc}]
\ref{thm:C2SBP-nec-cond-extinc-i} For any $n\in\N_0$,
\begin{align*}
E[Y_{n+1}|\bar{\mathcal{F}}_{n}]
&=E[E[\bar{F}_{n+1}|\bar{\mathcal{H}}_{n+1}]|\bar{\mathcal{F}}_{n}]-E[E[\bar{F}_{n+1}I_{\{\bar{M}_{n+1}=0\}}|\bar{\mathcal{H}}_{n+1}]|\bar{\mathcal{F}}_{n}]\\
&=E[\lambda\varphi_{n+1}(X_{n+1})|\bar{\mathcal{F}}_{n}]-E[\varphi_{n+1}(X_{n+1})\lambda^{\varphi_{n+1}(X_{n+1})}|\bar{\mathcal{F}}_{n}]\\
&=\lambda E[\epsilon(X_{n+1})|\bar{\mathcal{F}}_{n}]-\lambda E[h_{X_{n+1}}'(\lambda)|\bar{\mathcal{F}}_{n}]\quad a.s.,
\end{align*}
\noindent where $I_A$ stands for the indicator function of the set $A$.

\vspace{0.25cm}

\noindent \ref{thm:C2SBP-nec-cond-extinc-ii}-\ref{thm:C2SBP-nec-cond-extinc-ii-a} Taking into account \ref{thm:C2SBP-nec-cond-extinc-i} and the fact that the p.g.f. of the variable $\varphi_0(k)$ is $h_k(u)=(1-\rho+\rho u)^k$, we have that there exists a positive constant $C_1$ such that
\begin{align*}
E[Y_{n+1}|\bar{\mathcal{F}}_{n}]
&=\lambda\rho E[X_{n+1}|\bar{\mathcal{F}}_{n}]-\lambda E[E[h_{X_{n+1}}'(\lambda)|\bar{\mathcal{G}}_{n+1}]|\bar{\mathcal{F}}_{n}]\\
&=\lambda\rho\mu Y_n-\rho\lambda E\left[X_{n+1}(1-\rho+\rho\lambda)^{X_{n+1}-1}|\bar{\mathcal{F}}_{n}\right]\\
&=\lambda\rho\mu Y_n-\rho\lambda Y_n f(1-\rho+\rho\lambda)^{Y_n-1} f'(1-\rho+\rho\lambda)\\
&= \lambda\rho\mu Y_n- C_1 Y_n f(1-\rho+\rho\lambda)^{Y_n}\quad a.s.
\end{align*}
Analogously, we have
\begin{align*}
E[Y_{n+1}^2|\bar{\mathcal{F}}_{n}]&=E[E[\bar{F}_{n+1}^2|\bar{\mathcal{H}}_{n+1}]|\bar{\mathcal{F}}_{n}]-E[E[\bar{F}_{n+1}^2I_{\{\bar{M}_{n+1}=0\}}|\bar{\mathcal{H}}_{n+1}]|\bar{\mathcal{F}}_{n}]\\
&=E[\lambda(1-\lambda)\varphi_{n+1}(X_{n+1})|\bar{\mathcal{F}}_{n}]+E[\lambda^2\varphi_{n+1}(X_{n+1})^2|\bar{\mathcal{F}}_{n}]\\
&\phantom{=}-E[\varphi_{n+1}(X_{n+1})^2\lambda^{\varphi_{n+1}(X_{n+1})}\ |\bar{\mathcal{F}}_{n}]\\
&=\lambda(1-\lambda)\rho \mu Y_n +\lambda^2 \rho(1-\rho) \mu Y_n + \lambda^2\rho^2\delta^2 Y_n  +\lambda^2\rho^2 \mu^2Y_n^2-\lambda^2E[h_{X_{n+1}}''(\lambda) |\bar{\mathcal{F}}_{n}]\\
&\phantom{=}-\lambda E[h_{X_{n+1}}'(\lambda)|\bar{\mathcal{F}}_{n}]\\
&=\lambda(1-\lambda)\rho \mu Y_n +\lambda^2 \rho(1-\rho) \mu Y_n + \lambda^2\rho^2\delta^2 Y_n  +\lambda^2\rho^2 \mu^2 Y_n^2\\
&\phantom{=}-\lambda^2 \rho^2 Y_n(Y_n-1) f(1-\rho+\rho\lambda)^{Y_n-2}f'(1-\rho+\rho\lambda)^2\\
&\phantom{=} -\lambda^2\rho^2 Y_n f(1-\rho+\rho\lambda)^{Y_n-1}f''(1-\rho+\rho\lambda)\\
&\phantom{=}-\lambda \rho Y_n f(1-\rho+\rho\lambda)^{Y_n-1}f'(1-\rho+\rho\lambda).
\end{align*}

Combining this with the expression of the conditional expectation we have that, there exist positive constants $C_2$ and $C_3$ such that
\begin{align*}
\V[Y_{n+1}|\bar{\mathcal{F}}_{n}]
&\leq C_2 Y_n+ C_3 Y_n^2f(1-\rho+\rho\lambda)^{Y_n}
\end{align*}


In order to prove \ref{thm:C2SBP-nec-cond-extinc-ii}-\ref{thm:C2SBP-nec-cond-extinc-ii-b}, we can choose $0<\zeta<\lambda\mu\rho-1$, so that $1<\eta=\lambda\mu\rho-\zeta$.
%
%
Let us also define the  sets $A_n=\{\eta Y_n<Y_{n+1}\}$, $n\in\N_0$. Then, it is not difficult to verify that,
\begin{align*}
P(Y_n\to\infty|Y_0=y)&\geq \lim_{n\to\infty}P\left(\cap_{l=0}^{n}A_l|Y_0=y\right)\nonumber\\
&=\lim_{n\to\infty}P(A_0|Y_0=y)\prod_{l=1}^{n}P\left(A_l\big|\cap_{k=0}^{l-1}A_k\cap\{Y_0=y\}\right).
\end{align*}

Moreover, for each $l\in\N_0$ fixed, let us define the sets $B_{il}=\{Y_l=i\}\cap\big(\cap_{k=0}^{l-1}A_k\cap\{Y_0=y\}\big)$, then it is immediate that
$$\cup_{i=1}^\infty  B_{il}=\cap_{k=0}^{l-1}A_k\cap\{Y_0=y\},$$
and that if $\omega\in \cap_{k=0}^{l-1}A_k\cap\{Y_0=y\}$, then $Y_l(\omega)>\eta Y_{l-1}(\omega)>\ldots>\eta^ly$. Thus,
\begin{align*}
P\Big(A_l |\cap_{k=0}^{l-1}A_k\cap\{Y_0=y\}\Big)&=P\Big(A_l |\cup_{i=1}^\infty  B_{il}\Big)\geq \inf_{i>y\eta^l} P(A_l| Y_l=i)=\inf_{i>y\eta^l} P(A_0| Y_0=i),
\end{align*}
and now, we shall obtain a lower bound for $P(A_0| Y_0=i)$. Note that since, $f(1-\rho+\rho\lambda)^{i}\to 0$, as $i\to\infty$, then there exist $I_0\in\N$ and $0<\epsilon<\zeta$  such that $\rho\lambda f(1-\rho+\rho\lambda)^{i-1} f'(1-\rho+\rho\lambda)< \epsilon <\zeta$, for $i\geq I_0$.
As a consequence, on the one hand, by using Chebyshev's inequality, for $i\geq I_0$,
\begin{align*}\label{eq:prop:survival-non-iso4}
P\big(A_0^c|Y_0=i\big)
&=P\big(\zeta Y_0\leq \lambda\mu\rho Y_0- Y_{1} |Y_0=i\big)\nonumber\\
&\leq P\big((\zeta -\rho\lambda f(1-\rho+\rho\lambda)^{i-1} f'(1-\rho+\rho\lambda))Y_0 \leq |E[Y_{1}|Y_0=i]- Y_{1}|\big|Y_0=i\big)\nonumber\\
&\leq \frac{\V[Y_{1}|Y_0=i]}{(\zeta -\epsilon)^2 i^2}.
\end{align*}
On the other hand, using \ref{thm:C2SBP-nec-cond-extinc-ii-a} we have that there exist positive constants $K_1$ and $K_2$ satisfying
\begin{equation*}
P\big(A_0^c|Y_0=i\big)\leq \frac{K_1}{i}+K_2 f(1-\rho+\rho\lambda)^{i},\qquad \text{ for $i\geq I_0$}.
\end{equation*}

Note that since $\eta>1$, by taking $y\geq I_0$ we have $y\eta^l\geq I_0$, for all $l\in\N_0$. Thus, from all the above we obtain
\begin{align*}
P(Y_n\to\infty|Y_0=y)&\geq P(A_0|Y_0=y)\cdot\prod_{l=1}^{\infty}\bigg(1-\frac{K_1}{y\eta^l}-K_2 f(1-\rho+\rho\lambda)^{y\eta^l}\bigg)>0,\quad y\geq I_0.
\end{align*}

Finally, we observe that $\N$ is a class of communicating states for this process and hence, using the same arguments as in \cite{art-2004d}, p.47, we can prove that $P(Y_n\to\infty|Y_0=y)>0$ for each $y\in\N$.
\end{Prf}

\section{Proofs of the results in Section \ref{sec:BasicProperties}}

In the following results we make use of the next lemma. The proof is easy and it is omitted.

\begin{lema}\label{lem:females-males}
Let $n\in\N_0$ and $l_1,l_2\in\N_0$, then conditionally on $\{\phi_{n+1}(T_{n+1},\tT_{n+1})=l_1\}$, 
$$F_{n+1}\stackrel{d}{=}\sum_{i=1}^{l_1}f_{n i},\quad\text{ and }\quad M_{n+1}\stackrel{d}{=}l_1-\sum_{i=1}^{l_1}f_{n i}$$
where ``$\stackrel{d}{=}$'' denotes equal in distribution, and similarly,  conditionally on $\{\tphi_{n+1}(T_{n+1},\tT_{n+1})=l_2\}$, 
$$\tF_{n+1}\stackrel{d}{=}\sum_{i=1}^{l_2}\tilde{f}_{n i},\quad\text{ and }\quad \tM_{n+1}\stackrel{d}{=}l_2-\sum_{i=1}^{l_2}\tilde{f}_{n i},$$
where $\{f_{ni}: i\in\N, n\in\N_0\}$ and $\{\tilde{f}_{ni}: i\in\N, n\in\N_0\}$ are two independent families of r.v.s such that the former is a sequence of i.i.d.  following a Bernoulli distribution of parameter $\alpha$ and which are independent of $\phi_{n+1}(T_{n+1},\tT_{n+1})$ and the latter is a sequence of i.i.d. r.v.s following a Bernoulli distribution of parameter $\talpha$ and which are also independent of $\tphi_{n+1}(T_{n+1},\tT_{n+1})$. 
\end{lema}

\vspace{2ex}

\begin{Prf}[Proposition \ref{Moments FM|T}]
We only give the proofs for the predator females. The proofs are similar for the remaining conditional moments. 

\ref{Moments FM|T-i} Using the fact that $\mathcal{G}_{n}\subseteq \mathcal{H}_{n}$ and  Lemma \ref{lem:females-males}
\begin{align*}
E[F_{n}|\mathcal{G}_n]=E\left[E\left[\sum_{i=1}^{\phi_n(T_n,\tT_n)}f_{n-1 i}|\mathcal{H}_n\right]|\mathcal{G}_n\right]=\alpha E\left[\phi_n(T_n,\tT_n)|\mathcal{G}_n\right]=\alpha \varepsilon(T_n,\tT_n)\quad a.s.
\end{align*}

\vspace{-1ex}
\ref{Moments FM|T-ii} Taking into account that $\mathcal{G}_{n}\subseteq \mathcal{H}_{n}$ and Lemma \ref{lem:females-males}, we have that 
\begin{align*}
Var[F_{n}|\mathcal{G}_n]
&=Var\Bigg[E\Bigg[\sum_{i=1}^{\phi_n(T_n,\tilde{T}_n)}f_{ni}\Big|\mathcal{H}_{n}\Bigg]\Big|\mathcal{G}_{n}\Bigg]+E\Bigg[Var\Bigg[\sum_{i=1}^{\phi_n(T_n,\tilde{T}_n)}f_{ni}\Big|\mathcal{H}_{n}\Bigg]\Big|\mathcal{G}_{n}\Bigg]\\
&=Var[\alpha \phi_n(T_n,\tilde{T}_n)|\mathcal{G}_{n}]+E[\alpha(1-\alpha)\phi_n(T_n,\tilde{T}_n)|\mathcal{G}_{n}]\\
&=\alpha^2 \sigma^2(T_n,\tilde{T}_n)+\alpha(1-\alpha)\varepsilon(T_n,\tilde{T}_n) \quad a.s.
\end{align*}
\end{Prf}

\vspace{2ex}

\begin{Prf}[Proposition \ref{Moments FM|Z}]
We only provide the proof for the number of females in the predator population; the remaining  ones are similar.

\ref{Moments FM|Z-i} Taking into account $\mathcal{F}_{n}\subseteq \mathcal{G}_{n+1}$, Proposition \ref{Moments FM|T}~\ref{Moments FM|T-i}, the expectation of the control variables, and the independence of $T_n$ and $\tT_n$ given $Z_n$ and $\tZ_n$, we have that
\begin{align*}
E[F_{n+1}|\mathcal{F}_{n}]=E\left[\alpha\varepsilon(T_{n+1},\tT_{n+1})|\mathcal{F}_{n}\right]=\alpha E[T_{n+1}s(\tT_{n+1})|\mathcal{F}_{n}]=\alpha m Z_{n}E[s(\tT_{n+1})|\mathcal{F}_{n}].
\end{align*}

\vspace*{2ex}

\ref{Moments FM|Z-ii} By Proposition \ref{Moments FM|T}~\ref{Moments FM|T-i} and \ref{Moments FM|T-ii} and bearing in mind that $\mathcal{F}_{n}\subseteq \mathcal{G}_{n+1}$, we have that
\begin{align*}
Var[F_{n+1}|\mathcal{F}_{n}]&
=Var[\alpha T_{n+1}s(\tT_{n+1})|\mathcal{F}_{n}]+E[\alpha T_{n+1} s(\tT_{n+1})(1-\alpha s(\tT_{n+1}))|\mathcal{F}_{n}]\\
&=\alpha^2Var[ T_{n+1}s(\tT_{n+1})|\mathcal{F}_{n}]+\alpha E[T_{n+1} s(\tT_{n+1})|\mathcal{F}_{n}]-\alpha^2 E[T_{n+1}s(\tT_{n+1})^2|\mathcal{F}_{n}].
\end{align*}

On the one hand, for the second and third terms by using Proposition \ref{Moments Individuals}
$$E[T_{n+1}s(\tT_{n+1})|\mathcal{F}_{n}]=m Z_n E[s(\tT_{n+1})|\mathcal{F}_{n}],\quad \text{ and }\quad E[ T_{n+1}s(\tT_{n+1})^2|\mathcal{F}_{n}]=m Z_nE[s(\tT_{n+1})^2|\mathcal{F}_{n}].$$
Analogously, for the first term
\begin{align*}
Var[T_{n+1}s(\tT_{n+1})|\mathcal{F}_{n}]
&=Var[T_{n+1}|\mathcal{F}_{n}]E[s(\tT_{n+1})^2|\mathcal{F}_{n}]+E[T_{n+1}|\mathcal{F}_{n}]^2 Var[s(\tT_{n+1})|\mathcal{F}_{n}]\\
&=\sigma^2 Z_n E[s(\tT_{n+1})^2|\mathcal{F}_{n}]+m^2Z_n^2 Var[s(\tT_{n+1})|\mathcal{F}_{n}].
\end{align*}
Then, combining all the above the result follows.
\end{Prf}

\vspace{2ex}

\begin{Prf}[Proposition \ref{prop:expectation}]
By condition \ref{cond:A2}, Proposition \ref{Moments FM|Z} and the properties of the survival probability functions $s(\cdot)$ and $\ts(\cdot)$, we have that
\begin{align*}
E[Z_{n+1}|\mathcal{F}_n]
\leq E[F_{n+1}|\mathcal{F}_n]
=\alpha m Z_n E[s(\tT_{n+1})|\tZ_n]\leq \rho_2\alpha m Z_n\quad a.s.
\end{align*}
Analogously, we obtain that
$$E[\tZ_{n+1}|\mathcal{F}_n]\leq \tilde{\alpha}\tilde{m}\tZ_n E[\ts(T_{n+1})|\mathcal{F}_n]\leq \trho_2\tilde{\alpha}\tilde{m}\tZ_n\quad a.s.$$
\end{Prf}

\vspace{2ex}

\begin{Prf}[Proposition \ref{lem:moments-T-promiscuous}]
We provide the proofs for the number of predators; the arguments are similar for the number of preys.

\ref{lem:moments-T-promiscuous-i} Let us consider the $\sigma$-algebra $\mathcal{G}'_n=\sigma(Z_0,\tZ_0,T_l,\tT_l,Z_l,\tZ_l: l=1,\ldots,n)$, which satisfies $\mathcal{G}_n\subseteq \mathcal{G}'_n$, for $n\in\mathbb{N}$. Then, for $n\in\N$
\begin{align*}
E[T_{n+1}|\mathcal{G}_n]=E[E[T_{n+1}|\mathcal{G}'_n]|\mathcal{G}_n]
=m E[Z_n|\mathcal{G}_n]
=\alpha m T_ns(\tT_n)-\alpha m T_ns(\tT_n)(1-s(\tT_n)+\alpha s(\tT_n))^{T_n-1}.
\end{align*}
where we have used that 
\begin{align*}
E[Z_n|\mathcal{G}_n]
=E[F_n|\mathcal{G}_n]-E[F_nI_{\{M_n=0\}}|\mathcal{G}_n]
=\alpha T_ns(\tT_n)-\alpha h'_{(T_n,\tT_n)}(\alpha).
\end{align*}
with $h_{(t,\tilde{t})}(\cdot)$ being the p.g.f. of the variable $\phi_0(t,\tilde{t})$, i.e., $h_{(t,\tilde{t})}(u)=(1-s(\tilde{t})+s(\tilde{t})u)^t$, for $u\in [0,1]$, $t,\tilde{t}\in\N_0$.

\vspace{1ex}

\ref{lem:moments-T-promiscuous-ii} Since
$$E[T_{n+1}^2|\mathcal{G}'_n]=Var\left[\sum_{i=1}^{Z_n}t_{ni}\big|\mathcal{G}'_n\right]+E\left[\sum_{i=1}^{Z_n}t_{ni}\big|\mathcal{G}'_n\right]^2=Z_n\sigma^2+Z_n^2m^2,$$
we have
\begin{align*}
E[T_{n+1}^2|\mathcal{G}_n]= E[E[T_{n+1}^2|\mathcal{G}'_n]|\mathcal{G}_n]= E[Z_n\sigma^2+Z_n^2m^2|\mathcal{G}_n]= \sigma^2E[Z_n|\mathcal{G}_n]+m^2E[Z^2_n|\mathcal{G}_n].
\end{align*}

Now, bearing in mind that we computed $E[Z_n|\mathcal{G}_n]$ and $E[T_{n+1}|\mathcal{G}_n]$ in the proof of \ref{lem:moments-T-promiscuous-i}, it only remains to determine $E[Z^2_n|\mathcal{G}_n]$, and this is done as follows. First, note that $E[Z^2_n|\mathcal{G}_n]=E[F_n^2|\mathcal{G}_n]-E[F_n^2I_{\{M_n=0\}}|\mathcal{G}_n]$, and by Proposition \ref{Moments FM|T} we have
\begin{align*}
E[F_n^2|\mathcal{G}_n]&=Var[F_n|\mathcal{G}_n]+E[F_n|\mathcal{G}_n]^2=\alpha^2T_ns(\tT_n)(1-s(\tT_n))+\alpha(1-\alpha)T_n s(\tT_n)+\alpha^2 T_n^2 s(\tT_n)^2.
\end{align*}
Moreover,
\begin{align*}
E[F_n^2I_{\{M_n=0\}}|\mathcal{G}_n]
=\alpha^2 E[\phi_n(T_n,\tT_n)^2\alpha^{\phi_n(T_n,\tT_n)-2}|\mathcal{G}_n]
=\alpha^2 h''_{(T_n,\tT_n)}(\alpha)+\alpha h'_{(T_n,\tT_n)}(\alpha),
\end{align*}
and consequently
\begin{align*}
E[Z^2_n|\mathcal{G}_n]
&=\alpha^2T_ns(\tT_n)(1-s(\tT_n))+\alpha(1-\alpha)T_n s(\tT_n)+\alpha^2 T_n^2 s(\tT_n)^2\\
&\phantom{=}-\alpha^2 T_n(T_n-1)s(\tT_n)^2(1-s(\tT_n)+\alpha s(\tT_n))^{T_n-2}\\
&\phantom{=}-\alpha T_ns(\tT_n)(1-s(\tT_n)+\alpha s(\tT_n))^{T_n-1}
\end{align*}
Combining all the above the result yields.
\end{Prf}

\begin{Prf}[Proposition \ref{prop:extexp}]
Part \ref{prop:extexp-i} is easily deduced by applying well-known results about general Markov chains theory, bearing in mind that the states $(k,\tilde{k})$, with $(k,\tilde{k})\in\N_0^2\backslash\{(0,0)\}$ are transient (see, for example, \cite[Section I.17]{Chung}). 

\vspace*{0.25cm}

In order to proof \ref{prop:extexp-ii} we observe that 
\begin{align*}
\Omega&=\{Z_n\to \infty, \tZ_n\to \infty\}\cup \{Z_n\to 0, \tZ_n\to \infty\}\cup\{Z_n\to \infty, \tZ_n\to 0\} \cup \{Z_n\to 0, \tZ_n\to 0\}\\ 
&\phantom{=}\cup A_1\cup A_2,
\end{align*}
where $A_1=\{Z_n\nrightarrow 0, Z_n\nrightarrow \infty\}$, and $A_2=\{\tZ_n\nrightarrow 0, \tZ_n\nrightarrow \infty\}$.


We shall start by proving that $P(A_1)=0$. Given $\omega\in A_1$, then there exist $0<B<\infty$ such that for all $n_0$, there exists $n\geq n_0$ satisfying $0<Z_n(\omega)\leq B$, i.e., 
\begin{align*}
\omega &\in \cap_{n_0=1}^{\infty}\cup_{n=n_0}^{\infty}\{0<Z_n\leq B\},
\end{align*}
then,
\begin{align*}
A_1\subseteq\cup_{B=1}^\infty \limsup_{n\to\infty}\ \{0<Z_n\leq B\},
\end{align*}
and we conclude that $P(A_1)=0$ by applying \ref{prop:extexp-i}. Similar arguments lead to $P(A_2)=0$, and thus, the result follows.

\end{Prf}

\section{Proofs of the results in Section \ref{sec:Isolated}}\label{ape:isolated}

\begin{Prf}[Proposition \ref{prop:predator fixation}]
First, observe that 
\begin{align*}
P_{(i,j)}(Z_n\to \infty, \tZ_n\to 0)
&\leq\sum_{n=0}^\infty P_{(i,j)}\left(\{\tZ_n=0\} \cap (\cap_{k=0}^\infty Z_k>0\})\right)\\
&\leq
\sum_{n=0}^\infty P_{(i,j)}\left(Z_{n+1}>0|\tZ_{n}=0\right)P_{(i,j)}(\tZ_n=0).
\end{align*}
Now, taking into account the definition of the model and the fact that $\phi_0(t,0)=0$ a.s. for each $t\in\N_0$, by condition \ref{cond:A3} we conclude that for any $n\in\N_0$,
\begin{align*}
P_{(i,j)}(Z_{n+1}>0|\tZ_{n}=0)=P_{(i,j)}\Big(L(F_{n+1},\phi_{n+1}(T_{n+1},\tT_{n+1})-F_{n+1})>0|\tZ_{n}=0\Big)=0.
\end{align*}
\end{Prf}

\begin{Prf}[Proposition \ref{prop:impossible coexistence}]
We follow the same arguments as in the proof of Lemma 2 in \cite{art-bisexual-extinction-ylinked}. On the one hand, by Proposition \ref{prop:expectation} we have that $E[\tZ_{n+1}|\mathcal{F}_{n}]\leq \talpha \tilde{m} \tZ_{n}E[\ts(T_{n})|\mathcal{F}_n]$ a.s.

On the other hand, by the definition of the model and condition \ref{cond:A2}, $Z_n=L(F_n,M_n)\leq F_n \leq T_n$, for each $n\in\N$, and taking into account that the function $\ts(\cdot)$ is strictly decreasing we have that, for all $n\geq 1$, $\ts(Z_n)\geq \ts(T_n)$, and then
$$E[\tZ_{n+1}|\mathcal{F}_n]\leq \talpha \tilde{m} \tZ_{n}E[\ts(T_{n})|\mathcal{F}_n]\leq\talpha \tilde{m} \tZ_{n}\ts(Z_{n}).$$

Since $\lim_{x\to\infty}\ts(x)=0$, there exists $A>0$ such that $|x|\geq A$ implies $\ts(x)\leq \frac{1}{\talpha \tilde{m}}$, and then
\begin{equation}\label{equ:mart Z}
E[\tZ_{n+1}|\mathcal{F}_n]\leq \tZ_{n} \ \mbox{ a.s.\quad on } \{Z_n\geq A\}.
\end{equation}

We shall prove now that this implies that $P_{(i,j)}(Z_n\to\infty,\tZ_n\to\infty)=0$. To that end, we shall prove that, for every $N>0$,
\begin{equation}\label{equ:prob sup}
P_{(i,j)}\left(\left\{\inf_{n\geq N}Z_n > A\right\}\cap \{Z_n\to\infty,\tZ_n\to\infty\}\right)=0.
\end{equation}

Let us fix $N>0$ and define the stopping time:
$$T(A)=
\begin{cases}
\infty, & \mbox{ if } \inf_{n\geq N}Z_n \geq A,\\
\min\{n\geq N: Z_n < A\}, & \mbox{ otherwise},
\end{cases}
$$
and define also the sequence of r.v.s $\{Y_n\}_{n\in\N_0}$ as follows:
$$Y_n=\begin{cases}\tZ_{N+n}, & \mbox{ if } N+n\leq T(A),\\
\tZ_{T(A)}, & \mbox{ if } N+n>T(A).\end{cases}$$

To obtain \eqref{equ:prob sup}, we show that $\{Y_n\}_{n\in\N_0}$ is a non-negative super-martingale with respect to $\{\mathcal{F}_{N+n}\}_{n\in\N_0}$. Indeed, the variable $Y_n$ is $\mathcal{F}_{N+n}$-measurable for any $n\geq 0$. 

Let us fix $n\geq 0$, if $Z_{N+k}\geq A$, for each $k=0,\ldots,n$, then $T(A)\geq N+n+1$ and from \eqref{equ:mart Z} we obtain that
$$E[Y_{n+1}|\mathcal{F}_{N+n}]=E[\tZ_{N+n+1}|\mathcal{F}_{N+n}]\leq \tZ_{N+n}=Y_n \ \mbox{a.s.\quad on } \{Z_{N+k}\geq A: k=0,\ldots,n\}.$$

If there exists $k\in\{1,\ldots,n\}$ such that $Z_N\geq A,\ldots,Z_{N+k-1}\geq A$ and $Z_{N+k}<A$, then $T(A)\leq N+k<N+n+1$ and also
$$E[Y_{n+1}|\mathcal{F}_{N+n}]=E[\tZ_{T(A)}|\mathcal{F}_{N+n}]=Y_n \ \mbox{a.s.\quad on } \{Z_N\geq A,\ldots,Z_{N+k-1}\geq A,Z_{N+k}<A\}.$$

Finally, if $Z_N<A$, then $T(A)=N<N+n+1$ for all $n\geq 0$, and we get that
$$E[Y_{n+1}|\mathcal{F}_{N+n}]=E[\tZ_{N}|\mathcal{F}_{N+n}]=Y_n \ \mbox{a.s.\quad on } \{Z_{N}<A\}.$$

In short, since $B_n=\{\tZ_{N+k}\geq A: k=0,\ldots,n\} \in \mathcal{F}_{N+n}$ one deduces that
\begin{align*}
E[Y_{n+1}|\mathcal{F}_{N+n}]=E[Y_{n+1}|\mathcal{F}_{N+n}]I_{B_n}+E[Y_{n+1}|\mathcal{F}_{N+n}]I_{B_n^c}\leq \tZ_{N+n}I_{B_n} + \tZ_{T(A)}I_{B_n^c}=Y_n \mbox{ a.s.}
\end{align*}

Applying the martingale convergence theorem, we obtain the almost sure convergence of the sequence $\{Y_n\}_{n\in\N_0}$ to a non-negative and finite limit, $Y_\infty$ where
$$Y_{\infty}=\begin{cases}\lim_{n\to\infty}\tZ_n, & \mbox{ if } \inf_{n\geq N} Z_n >A,\\ \tZ_{T(A)}, & \mbox{ otherwise.}\end{cases}$$
Therefore we deduce \eqref{equ:prob sup} and the proof is completed.
\end{Prf}

\section{Proofs of the results in Section \ref{sec:Non-isolated}}\label{ape:non-isolated}

\begin{Prf}[Proposition \ref{prop:non-isol extinction}]
From Proposition \ref{prop:expectation} it is immediate to see that if $\rho_2m\alpha\leq 1$ and $\trho_2\tilde{m}\tilde{\alpha}\leq 1$, then $\{Z_n\}_{n\in\N_0}$ and $\{\tZ_n\}_{n\in\N_0}$ are both non-negative supermartingales with respect to the $\sigma$-algebra $\mathcal{F}_n$ and therefore both of them  converge a.s. to a finite limit. Now, the result is derived by Proposition \ref{prop:extexp}~\ref{prop:extexp-ii}. 
\end{Prf}

\vspace{0.5cm}

\begin{Prf}[Proposition \ref{prop:non-isol common extinction}]
On the one hand, if $\rho_1\alpha m>1$, then $P_{(i,j)}(Z_n\to\infty,\tZ_n\to 0)>0$  by Proposition \ref{prop: non isolated prey and predator fixation}~\ref{prop:non-iso-predator-fixation} and together with Proposition \ref{prop:extexp}~\ref{prop:extexp-ii} we obtain $P_{(i,j)}(Z_n\to 0,\tZ_n\to 0)<1$. The proof is analogous if $\trho_2\talpha \tilde{m}>1$ by Propositions \ref{prop: non isolated prey and predator fixation}~\ref{prop:non-iso-prey-fixation} and \ref{prop:extexp}~\ref{prop:extexp-ii}.

On the other hand, if $\max\{\rho_1\alpha m,\trho_2\talpha \tilde{m}\}\leq 1$, then $P_{(i,j)}(Z_n\to \infty,\tZ_n\to 0)=0$ by Proposition \ref{prop: non isolated prey and predator fixation}~\ref{prop:non-iso-predator-fixation} and $P_{(i,j)}(\tZ_n\to \infty)=0$ by Proposition~\ref{prop:non-isol extinction}~\ref{prop:non-isol extinction-ii}. The result yields by \eqref{equ:Ztinfty} and Proposition~\ref{prop:extexp}~\ref{prop:extexp-ii}.
\end{Prf}

\vspace{0.5cm}

\begin{Prf}[Theorem \ref{prop:non-isol coexistence}]
\ref{prop:non-isol coexistence-i} By Proposition \ref{prop:expectation}, we have that
\begin{equation*}
E[\tZ_{n+1}|\mathcal{F}_n]\leq \talpha \tm \tZ_n E\left[\ts(T_{n+1})|\mathcal{F}_n\right].
\end{equation*}
Taking into account that $\lim_{x\to\infty}\ts(x)=\trho_1$, we have that for $\nu=\frac{1-\trho_1\tm\talpha}{\tm\talpha}>0$, there exists $A>0$ such that $|x|\geq A$ implies $\ts(x)\leq \trho_1+\nu$, and then
\begin{equation*}\label{equ:mart-Z-non-iso}
E[\tZ_{n+1}|\mathcal{F}_n]\leq \tm \talpha \tZ_n (\trho_1+\nu)\leq \tZ_n\ \mbox{ a.s.\quad on } \{Z_n\geq A\}.
\end{equation*}
The result follows with the same arguments as in the proof of Proposition \ref{prop:impossible coexistence}.
%
\vspace{0.25cm}

\ref{prop:non-isol coexistence-ii} First, notice that $\{T_n\to\infty,\tT_n\to\infty\}=\{Z_n\to\infty,\tZ_n\to\infty\}$ a.s. and thus, it is enough to prove that $P_{(i,j)}(T_n\to\infty,\tT_n\to\infty)>0$ if $\min\{\rho_2 m\alpha,\trho_1\tm\talpha\}>1$. To do that, first we prove that there exist $I_1,J_1\in\N_0$ large enough so that $P_{(i,j)}(T_n\to\infty,\tT_n\to\infty)>0$, for each $i,j\in\N$ satisfying $i\geq I_1$ and $j\geq J_1$, and then we prove that this implies the desired result.

Due to the fact that $1<\alpha\rho_2 m$ and $1<\talpha\trho_1 \tilde{m}$, we can take $0<\varepsilon_1$ and $0<\varepsilon_2$ small enough so that $\eta_1=\alpha\rho_2 m-\varepsilon_1>1$ and $\eta_2=\talpha \trho_1\tilde{m}-\varepsilon_2>1$. On the other hand, since $p_0+p_1+p_2<1$ and $\tilde{p}_0+\tilde{p}_1+\tilde{p}_2<1$, there exist $k_1,k_2\geq 3$ such that $p_{k_1}>0$ and $\tilde{p}_{k_2}>0$ and then, $\eta_0=k_1-1>1$ and $\tilde{\eta}_0=k_2-1>1$. Now, let us define the following sets:
$$A_0=\{\eta_0Z_0<T_1, \tilde{\eta}_0\tZ_0<\tT_1\}\quad \mbox{ and }\quad A_n=\{\eta_1T_n<T_{n+1},\eta_2\tT_n<\tT_{n+1}\},\quad n\in\N.$$
Then, it is not difficult to verify that,
\begin{align*}\label{equ: non-aisolated coexistence probAn1}
P_{(i,j)}(T_n\to\infty,\tT_n\to\infty)&\geq P_{(i,j)}\left(\cap_{n=0}^{\infty}A_n\right)\\
&=\lim_{n\to\infty}P_{(i,j)}(A_0)\prod_{l=1}^{n}P\left(A_l|\cap_{k=0}^{l-1}A_k\cap \{Z_0=i,\tZ_0=j\}\right),
\end{align*}
and for each $l\in\N$ fixed we have that
\begin{align*}
P\big(A_l|\cap_{k=0}^{l-1}A_k&\cap \{Z_0=i,\tZ_0=j\}\big)=\\
&=P\left(A_{l}|\cup_{i_0=1}^\infty\cup_{j_0=1}^\infty\big(\cap_{k=0}^{l-1}A_k\cap \{(T_l,\tT_l)=(i_0,j_0)\}\cap \{Z_0=i,\tZ_0=j\}\big)\right)\nonumber\\
&\geq \inf_{i\eta_0\eta_1^{l-1}<i_0, \atop j\tilde{\eta}_0\eta_2^{l-1}<j_0}P\left(A_{l}|(T_l,\tT_l)=(i_0,j_0)\right)\nonumber\\
&=\inf_{i\eta_0\eta_1^{l-1}<i_0, \atop j\tilde{\eta}_0\eta_2^{l-1}<j_0}P(A_1|(T_1,\tT_1)=(i_0,j_0)).
\end{align*}

As a consequence, we only need to obtain convenient lower bounds for $P_{(i,j)}(A_0)$ and \linebreak$P(A_1|(T_1,\tT_1)=(i_0,j_0))$. For the former we have
\begin{align*}
P_{(i,j)}(T_1>\eta_0 Z_0,\tT_1>\tilde{\eta}_0 \tZ_0)
&\geq P(T_1=k_1 i|Z_0=i)\cdot P(\tT_1=k_2 j|\tZ_0=j)
\geq p_{k_1}^i\tilde{p}_{k_2}^j>0.
\end{align*}
%
%
%
%
%
%

For the latter, 
we shall obtain a convenient upper bound for $P(A^c_1|(T_1,\tT_1)=(i_0,j_0))$. 
First, from Proposition \ref{lem:moments-T-promiscuous} we have
\begin{align*}
P\big(T_2\leq \eta_1 T_1|(T_1,\tT_1)=(i_0,j_0)\big)
&=P\Big(i_0\big(\varepsilon_1-\alpha m \rho_2+\alpha ms(j_0)-\alpha m s(j_0)(1-s(j_0)+\alpha s(j_0))^{i_0-1}\big)\\
&\phantom{=}\hspace{2em}\leq |E[T_2|(T_1,\tT_1)=(i_0,j_0)]-T_2| \big|(T_1,\tT_1)=(i_0,j_0)\Big).
\end{align*}
Second, note that since $\lim_{x\to\infty} s(x)=\rho_2$, and $\lim_{x\to\infty} \sup_{j_0\in\N}\{(1-s(j_0)+\alpha s(j_0))^{x}\}=0$, we can choose $I_0\in\N$ and $J_0\in\N$ such that $|\alpha m\rho_2-\alpha m s(j_0)|<\varepsilon_1/4$, and $|\alpha m s(j_0)(1-s(j_0)+\alpha s(j_0))^{i_0-1}|<\varepsilon_1/4$, for any $i_0\geq I_0$ and $j_0\geq J_0$. Now, for $i_0\geq I_0$ and $j_0\geq J_0$, using Chebyschev's inequality we have
\begin{align*}
P\big(T_2\leq \eta_1 T_1|(T_1,\tT_1)=(i_0,j_0)\big)&\leq P\Big(\frac{i_0 \varepsilon_1}{2}\leq |E[T_2|(T_1,\tT_1)=(i_0,j_0)]-T_2| \big|(T_1,\tT_1)=(i_0,j_0)\Big)\\
&\leq \frac{4 Var[T_2|(T_1,\tT_1)=(i_0,j_0)]}{\varepsilon_1^2 i_0^2},
\end{align*}
and then by Proposition \ref{lem:moments-T-promiscuous} we have that there exist constants $K_1>0$ and $K_2>0$ satisfying
%
%
%
\begin{align*}
P\big(T_2\leq \eta_1 T_1|(T_1,\tT_1)=(i_0,j_0)\big)\leq \frac{K_1}{i_0}+K_2(1-(1-\alpha)\rho_1)^{i_0},\quad i_0\geq I_0,j_0\geq J_0.
\end{align*}

Analogously, we can prove that there exist $I_1\geq I_0$, $J_1\geq J_0$ and constants $K_3>0$ and $K_4>0$ satisfying that for $i_0\geq I_1$ and $j_0\geq J_1$, 
\begin{align*}
P\big(\tT_2\leq \eta_2 \tT_1|(T_1,\tT_1)=(i_0,j_0)\big)&\leq P\Big(\frac{j_0 \varepsilon_2}{2}\leq |E[\tT_2|(T_1,\tT_1)=(i_0,j_0)]-T_2| \big|(T_1,\tT_1)=(i_0,j_0)\Big)\\
&\leq \frac{K_3}{j_0}+K_4(1-(1-\talpha)\trho_1)^{j_0}.
\end{align*}

Combining all the above, if $i\geq I_1$ and $j\geq J_1$, then since $\eta_0,\tilde{\eta}_0,\eta_1,\eta_2>1$, we have $i\eta_0\eta_1^{l-1}\geq I_1$ and $j\tilde{\eta}_0\eta_2^{l-1}\geq J_1$, for all $l\in\N_0$, then it is immediate to verify that 
\begin{align*}
P_{(i,j)}(T_n\to\infty,\tT_n\to\infty)& \geq P_{(i,j)}(A_0)\cdot\prod_{l=1}^\infty\bigg(1-\frac{K_1}{i\eta_0\eta_1^{l-1}}-K_2(1-(1-\alpha)\rho_1)^{i\eta_0\eta_1^{l-1}}-\frac{K_3}{j\tilde{\eta}_0\eta_2^{l-1}}\\
&\phantom{\geq}-K_4(1-(1-\talpha)\trho_1)^{j\tilde{\eta}_0\eta_2^{l-1}}\bigg)>0,
\end{align*}
and then $P_{(i,j)}(Z_n\to\infty,\tZ_n\to\infty)>0$ for $i\geq I_1$ and $j\geq J_1$.
Finally, we prove that the result also holds for any $i,j\in\N$. Let us fix $i,j\in\N$ such that either $i<I_1$ or $j<J_1$. By Proposition \ref{prop:states}~\ref{prop:states-promiscuous-iii}, there exists $k_0\in\N$ such that $P(Z_{k_0}=I_1,\tZ_{k_0}=J_1|Z_0=i,\tZ_0=j)>0$. Then, we have
\begin{align*}
P(Z_n\to\infty,\tZ_n\to\infty|Z_0=i,\tZ_0=j)&\geq P(Z_n\to\infty,\tZ_n\to\infty|Z_{k_0}=I_1,\tZ_{k_0}=J_1)\\
&\phantom{=}\cdot P(Z_{k_0}=I_1,\tZ_{k_0}=J_1|Z_0=i,\tZ_0=j)>0.
\end{align*}

\end{Prf}


\end{document}